\documentclass[11pt,a4paper]{amsart}
\usepackage[all]{xy}

\newcommand{\myref}[1]{(\ref{#1})}

\newtheorem{theorem}{Theorem}[section]
\newtheorem{corollary}[theorem]{Corollary}
\newtheorem{lemma}[theorem]{Lemma}
\newtheorem{proposition}[theorem]{Proposition}
\newtheorem{conjecture}[theorem]{Conjecture}

\theoremstyle{definition}

\newcommand{\ie}{{\em i.e.}\ }
\newcommand{\cf}{{\em cf.}\ }
\newcommand{\eg}{{\em e.g.}\ }
\newcommand{\ko}{\: , \;}
\newcommand{\ul}[1]{\underline{#1}}
\newcommand{\ol}[1]{\overline{#1}}
\renewcommand{\tilde}[1]{\widetilde{#1}}

\newcommand{\Z}{\mathbb{Z}}
\newcommand{\N}{\mathbb{N}}
\renewcommand{\S}{\mathbb{S}}
\newcommand{\Q}{\mathbb{Q}}

\newcommand{\ca}{{\mathcal A}}
\newcommand{\cb}{{\mathcal B}}
\newcommand{\cc}{{\mathcal C}}
\newcommand{\cd}{{\mathcal D}}
\newcommand{\ce}{{\mathcal E}}

\newcommand{\cg}{{\mathcal G}}
\newcommand{\ch}{{\mathcal H}}

\newcommand{\cl}{{\mathcal L}}
\newcommand{\cm}{{\mathcal M}}
\newcommand{\cn}{{\mathcal N}}
\newcommand{\co}{{\mathcal O}}

\newcommand{\cR}{{\mathcal R}}

\newcommand{\cs}{{\mathcal S}}
\newcommand{\ct}{{\mathcal T}}
\newcommand{\cu}{{\mathcal U}}

\newcommand{\eps}{\varepsilon}

\newcommand{\opname}[1]{\operatorname{#1}}

\newcommand{\ra}{\rightarrow}
\newcommand{\la}{\leftarrow}
\newcommand{\iso}{\stackrel{_\sim}{\rightarrow}}
\newcommand{\id}{\mathbf{1}}

\newcommand{\op}{^{op}}
\newcommand{\ten}{\otimes}
\newcommand{\lten}{\overset{L}{\ten}}
\newcommand{\obj}{\operatorname{obj}}
\newcommand{\HOM}{{\mathcal H}om}
\newcommand{\RHOM}{{\mathcal R}{\mathcal H}om}
\newcommand{\Ext}{\opname{Ext}}
\newcommand{\Hom}{\opname{Hom}}
\newcommand{\Aut}{\opname{Aut}}
\newcommand{\Map}{\opname{Map}}
\newcommand{\rep}{\opname{rep}}
\newcommand{\Ho}{\opname{Ho}}
\newcommand{\Sset}{\opname{Sset}}

\newcommand{\bp}{{\mathbf p}}
\newcommand{\bi}{{\mathbf i}}
\newcommand{\R}{{\mathbf R}}
\renewcommand{\L}{{\mathbf L}}

\newcommand{\Mod}{\opname{Mod}\nolimits}
\renewcommand{\mod}{\opname{mod}\nolimits}
\newcommand{\dgcat}{\opname{dgcat}}
\newcommand{\per}{\opname{per}}
\newcommand{\parf}{\opname{par}}
\newcommand{\tria}{\opname{tria}}
\newcommand{\Hqe}{\opname{Hqe}}
\newcommand{\Hmo}{\opname{Hmo}}
\newcommand{\ac}{{\mathcal A}c}

\newcommand{\Spec}{\opname{Spec}}
\newcommand{\pretr}{\opname{pretr}}
\newcommand{\colim}{\opname{colim}}
\newcommand{\hocolim}{\opname{hocolim}}

\newcommand{\derh}{\cd\ch}
\newcommand{\HHs}{{HH}^*}

\title[On dg categories]{On differential graded categories}

\author{Bernhard Keller}

\address{UFR de Math\'ematiques, UMR 7586 du CNRS, Case 7012,
   Universit\'e Paris 7, 2 place Jussieu, 75251 Paris Cedex 05,
   France}
\email{keller@math.jussieu.fr}

\begin{document}

\begin{abstract}
Differential graded categories enhance our understanding
of triangulated categories appearing in algebra and geometry.
In this survey, we review their foundations and report on recent
work by Drinfeld, Dugger-Shipley, \ldots, To\"en and To\"en-Vaqui\'e.
\end{abstract}


 \date{February 26, 2006}


\maketitle


\section{Introduction}

\subsection{Triangulated categories and dg categories}
Derived categories were invented by Grothendieck-Verdier in the
early sixties in order to formulate Grothendieck's duality theory
for schemes, \cf \cite{Illusie90}. Today, they have become an
important tool in many branches of algebraic geometry, in
algebraic analysis, non commutative algebraic geometry,
representation theory, mathematical physics \ldots\ . In an
attempt to axiomatize the properties of derived categories,
Grothendieck-Verdier introduced the notion of a triangulated
category. For a long time, triangulated categories were considered
too poor to allow the development of more than a rudimentary
theory. This vision has changed in recent years \cite{Neeman99}
\cite{Neeman05}, but the fact remains that many important
constructions of derived categories cannot be performed with
triangulated categories. Notably, tensor products and functor
categories formed from triangulated categories are no longer
triangulated. One approach to overcome these problems has been the
theory of derivators initiated by Heller \cite{Heller88} and
Grothendieck \cite{Grothendieck90}, \cf also \cite{Keller91}, at
the beginning of the nineties. Another, perhaps less formidable
one is the theory of differential graded categories (=dg
categories), together with its cousin, the theory of
$A_\infty$-categories.

Dg categories already appear in \cite{Kelly65}. In the seventies,
they found applications \cite{Roiter78} \cite{Drozd80} in the
representation theory of finite-dimensional algebras. The idea to
use dg categories to `enhance' triangulated categories goes back
at least to Bondal-Kapranov \cite{BondalKapranov90}, who were
motivated by the study of exceptional collections of coherent
sheaves on projective varieties.

The synthesis of Koszul duality
\cite{BeilinsonGinzburgSchechtman88}
\cite{BeilinsonGinzburgSoergel96} with Morita theory for derived
categories \cite{Rickard89} was the aim of the study of the
unbounded derived category of a dg category in \cite{Keller94}.

It is now well-established that invariants like $K$-theory,
Hochschild (co-)homolo\-gy and cyclic homology associated with a
ring or a variety `only depend' on its derived category. However,
in most cases, the derived category (even with its triangulated
structure) is not enough to compute the invariant, and the datum
of a triangle equivalence between derived categories is not enough
to construct an isomorphism between invariants (\cf
Dugger-Shipley's \cite{DuggerShipley06a} results in
section~\ref{ss:topological}). Differential graded categories
provide the necessary structure to fill this gap. This idea was
applied to $K$-theory by Thomason-Trobaugh
\cite{ThomasonTrobaugh90} and to cyclic homology in
\cite{Keller98} \cite{Keller99}.

The most useful operation which {\em can be performed} on
triangulated categories is the passage to a Verdier quotient. It
was therefore important to lift this operation to the world of
differential graded categories. This was done implicitly in
\cite{Keller99} but explicitly, by Drinfeld, in \cite{Drinfeld04}.

In a certain sense, differential graded categories and
differential graded functors contain too much information and the
main problem in working with them consists in `discarding what is
irrelevant'. It now appears clearly that the best tool for doing
this are Quillen model categories \cite{Quillen67}: They provide a
homotopy theoretic framework which allows simple, yet precise
statements and rigorous but readable proofs. Building on the
techniques of \cite{Drinfeld04}, a suitable model structure on the
category of small differential graded categories was constructed
in \cite{Tabuada05a}. Starting from this structure,  To\"en has
given a new approach to Morita theory for dg categories
\cite{Toen04}. In their joint work \cite{ToenVaquie05}, To\"en and
Vaqui\'e have applied this to the construction of moduli stacks of
objects in dg categories, and notably in categories of perfect
complexes arising in geometry and representation theory.

Thanks to \cite{Toen04}, \cite{KockToen05} and to recent work by
Tamarkin \cite{Tamarkin05}, we are perhaps getting closer to
answering Drinfeld's question \cite{Drinfeld04}: {\em What do DG
categories form?}

\subsection{Contents} After introducing notations and basic
definitions in section~\ref{s:definition} we review the derived
category of a dg category in section~\ref{s:derived-category}.
This is the first opportunity to practice the language of model
categories. We present the structure theorems for algebraic
triangulated categories which are compactly generated or, more
generally, well-generated. We conclude with a survey of recent
important work by Dugger and Shipley on topological Morita
equivalence for dg categories. In
section~\ref{s:HtpyCatSmalldgCat}, we present the homotopy
categories of dg categories and of `triangulated' dg categories
following To\"en's work \cite{Toen04}. The most important points
are the description of the mapping spaces of the homotopy category
via quasi-functors (Theorem~\ref{thm:mapping-spaces}), the closed
monoidal structure (Theorem~\ref{thm:closed-monoidal}) and the
characterization of dg categories of finite type
(Theorem~\ref{thm:finite-type}). We conclude with a summary of the
applications to moduli problems. In the final
section~\ref{s:invariants}, we present the most important
invariance results for $K$-theory, Hochschild (co-)homology and
cyclic homology. The derived Hall algebra presented in
section~\ref{ss:Hall-algebras} is a new invariant due to To\"en
\cite{Toen05}. Its further development might lead to significant
applications in representation theory.


\subsection{Acknowledgments} I thank Bertrand To\"en, Henning Krause,
Brooke Shipley and Gon\c{c}alo Tabuada for helpful comments on
previous versions of this article.

\section{Definition}
\label{s:definition}

\subsection{Notations} \label{ss:notations}
Let $k$ be a commutative ring, for example a field or the ring of
integers $\Z$. We will write $\otimes$ for the tensor product over
$k$. Recall that a {\em $k$-algebra} is a $k$-module $A$ endowed
with a $k$-linear associative multiplication $A\otimes_k A \to A$
admitting a two-sided unit $1\in A$. For example, a $\Z$-algebra
is just a (possibly non commutative) ring. A {\em $k$-category
$\ca$} is a `$k$-algebra with several objects' in the sense of
Mitchell \cite{Mitchell72}. Thus, it is the datum of a class of
objects $\obj(\ca)$, of a $k$-module $\ca(X,Y)$ for all objects
$X$, $Y$ of $\ca$, and of $k$-linear associative composition maps
\[
\ca(Y,Z)\otimes \ca(X,Y) \to \ca(X,Z) \ko (f,g) \mapsto fg
\]
admitting units $\id_X\in\ca(X,X)$. For example, we can view
$k$-algebras as $k$-categories with one object. The category $\Mod
A$ of right $A$-modules over a $k$-algebra $A$ is an example of a
$k$-category with many objects. It is also an example of a {\em
$k$-linear} category, \ie a $k$-category which admits all finite
direct sums.

A {\em graded $k$-module} is a $k$-module $V$ together with a
decomposition indexed by the positive and the negative integers:
\[
V=\bigoplus_{p\in\Z} V^p.
\]
The {\em shifted module $V[1]$} is defined by $V[1]^p=V^{p+1}$,
$p\in\Z$. A {\em morphism} $f: V \to V'$ of graded $k$-modules of
degree $n$ is a $k$-linear morphism such that $f(V^p)\subset V^{p+n}$
for all $p\in\Z$. The {\em tensor product $V\ten W$} of two graded
$k$-modules $V$ and $W$ is the graded $k$-module with components
\[
(V\ten W)^n = \bigoplus_{p+q=n} V^p \ten W^q \ko n\in\Z.
\]
The {\em tensor product $f\ten g$} of two maps $f: V \to V'$ and
$g: W \to W'$ of graded $k$-modules is defined using the {\em
Koszul sign rule}: We have
\[
(f\ten g)(v\ten w) = (-1)^{pq} f(v)\ten g(w)
\]
if $g$ is of degree $p$ and $v$ belongs to $V^q$. A {\em graded
$k$-algebra} is a graded $k$-module $A$ endowed with a
multiplication morphism $A\ten A \to A$ which is graded of degree
$0$, associative and admits a unit $1\in A^0$. We identify
`ordinary' $k$-algebras with graded $k$-algebras concentrated in
degree $0$. We write $\cg(k)$ for the {\em category of graded
$k$-modules}.

A {\em differential graded (=dg) $k$-module} is a $\Z$-graded
$k$-module $V$ endowed with a {\em differential $d_V$}, \ie a map
$d_V: V \to V$ of degree $1$ such that $d_V^2=0$. Equivalently,
$V$ is a {\em complex} of $k$-modules. The {\em shifted dg module
$V[1]$} is the shifted graded module endowed with the differential
$-d_V$. The {\em tensor product} of two dg $k$-modules is the
graded module $V\ten W$ endowed with the differential $d_V\ten
\id_W + \id_V \ten d_W$.

\subsection{Differential graded categories} \label{ss:diffgradedcat}
A {\em differential graded} or {\em dg category} is a $k$-category
$\ca$ whose morphism spaces are dg $k$-modules and whose
compositions
\[
\ca(Y,Z)\ten \ca(X,Y) \to \ca(X,Z)
\]
are morphisms of dg $k$-modules.

For example, dg categories with one object may be identified with
{\em dg algebras}, \ie graded $k$-algebras endowed with a
differential $d$ such that the Leibniz rule holds:
\[
d(fg) = d(f) g + (-1)^p f\,d(g)
\]
for all $f\in A^p$ and all $g$. In particular, each ordinary
$k$-algebra yields a dg category with one object. A typical
example with several objects is obtained as follows: Let $B$ be a
$k$-algebra and $\cc(B)$ the category of complexes of right
$B$-modules
\[
\xymatrix{ \ldots \ar[r] & M^p \ar[r]^{d_M} & M^{p+1} \ar[r] &
\ldots \ko p\in\Z}.
\]
For two complexes $L,M$ and an integer $n\in\Z$, we define
$\HOM(L,M)^n$ to be the $k$-module formed by the morphisms $f: L
\to M$ of graded objects of degree $n$, \ie the families $f=(f^p)$
of morphisms $f^p : L^p \to M^{p+n}$, $p\in\Z$, of $B$-modules. We
define $\HOM(L,M)$ to be the graded $k$-module with components
$\HOM(L,M)^n$ and whose differential is the commutator
\[
d(f) = d_M \circ f - (-1)^n f\circ d_L.
\]
The {\em dg category $\cc_{dg}(B)$} has as objects all complexes
and its morphisms are defined by
\[
\cc_{dg}(B)(L,M)= \HOM(L,M).
\]
The composition is the composition of graded maps.

Let $\ca$ be a dg category. The {\em opposite dg category}
$\ca\op$ has the same objects as $\ca$ and its morphisms are
defined by
\[
\ca\op(X,Y)=\ca(Y,X);
\]
the composition of $f\in\ca\op(Y,X)^p$ with $g\in\ca\op(Z,Y)^q$ is
given by $(-1)^{pq} gf$. The {\em category $Z^0(\ca)$} has the
same objects as $\ca$ and its morphisms are defined by
\[
(Z^0 \ca)(X,Y) = Z^0(\ca(X,Y)) \ko
\]
where $Z^0$ is the kernel of $d: \ca(X,Y)^0 \to \ca(X,Y)^1$. The
{\em category $H^0(\ca)$} has the same objects as $\ca$ and its
morphisms are defined by
\[
(H^0(\ca))(X,Y) = H^0(\ca(X,Y)) \ko
\]
where $H^0$ denotes the $0$th homology of the complex $\ca(X,Y)$.
For example, if $B$ is a $k$-algebra, we have an isomorphism of
categories
\[
Z^0(\cc_{dg}(B)) = \cc(B)
\]
and an isomorphism of categories
\[
H^0(\cc_{dg}(B)) = \ch(B) \ko
\]
where $\ch(B)$ denotes the {\em category of complexes up to
homotopy}, \ie the category whose objects are the complexes and
whose morphisms are the morphisms of complexes modulo the
morphisms $f$ homotopic to zero, \ie such that $f=d(g)$ for some
$g\in \HOM(L,M)^{-1}$. The {\em homology category} $H^*(\ca)$ is
the graded category with the same objects as $\ca$ and
morphisms spaces $H^*\ca(X,Y)$.

\subsection{The category of dg categories} \label{ss:dgcat}
Let $\ca$ and $\ca'$ be dg categories. A {\em dg functor
$F:\ca\to\ca'$} is given by a map $F:\obj(\ca)\to\obj(\ca')$ and
by morphisms of dg $k$-modules
\[
F(X,Y) : \ca(X,Y) \to \ca(FX,FY) \ko X,Y\in\obj(\ca)\ko
\]
compatible with the composition and the units. The {\em category
of small dg categories $\dgcat_k$} has the small dg categories as
objects and the dg functors as morphisms. Note that it has an
initial object, the empty dg category $\emptyset$, and a final
object, the dg category with one object whose endomorphism ring is
the zero ring. The {\em tensor product $\ca\ten\cb$} of two dg
categories has the class of objects $\obj(\ca)\times\obj(\cb)$ and
the morphism spaces
\[
(\ca\ten\cb)((X,Y),(X',Y'))= \ca(X,X')\ten \cb(Y,Y')
\]
with the natural compositions and units.

For two dg functors $F,G: \ca \to \cb$, the {\em complex
of graded morphisms $\HOM(F,G)$} has as its $n$th component the
module formed by the families of morphisms
\[
\phi_X \in \cb(FX,GX)^n
\]
such that $(Gf)(\phi_X)=(\phi_Y)(Ff)$ for all $f\in\ca(X,Y)$,
$X,Y\in\ca$. The differential is induced by that of $\cb(FX,GX)$.
The set of {\em morphisms $F \to G$} is by definition in bijection
with $Z^0\HOM(F,G)$.

Endowed with the tensor product, the category $\dgcat_k$ becomes a
symmetric tensor category which admits an internal $\Hom$-functor,
\ie
\[
\Hom(\ca\ten\cb, \cc) = \Hom(\ca,\HOM(\cb,\cc)) \ko
\]
for $\ca,\cb,\cc\in\dgcat_k$, where $\HOM(\cb,\cc)$ has the dg
functors as objects and the morphism space $\HOM(F,G)$ for two dg
functors $F$ and $G$. The unit object is the dg category
associated with the $k$-algebra $k$.

A {\em quasi-equivalence} is a dg functor $F: \ca\to\ca'$ such
that
\begin{itemize}
\item[1)] $F(X,Y)$ is a quasi-isomorphism for all objects $X$, $Y$ of $\ca$
and
\item[2)] the induced functor $H^0(F): H^0(\ca)\to H^0(\ca')$ is an
equivalence.
\end{itemize}
Note that neither the tensor product nor the internal
$\Hom$-functor respect the quasi-equivalences, a source of
technical difficulties.

\section{The derived category of a dg category}
\label{s:derived-category}

\subsection{Dg modules} \label{ss:dgmodules}
Let $\ca$ be a small dg category. A {\em left dg $\ca$-module} is
a dg functor
\[
L : \ca \to \cc_{dg}(k)
\]
and a {\em right dg $\ca$-module} a dg functor
\[
M : \ca\op\to \cc_{dg}(k).
\]
Equivalently, a right dg $\ca$-module $M$ is given by complexes
$M(X)$ of $k$-modules, for each $X\in\obj(\ca)$, and by morphisms
of complexes
\[
M(Y) \ten \ca(X,Y) \to M(X)
\]
compatible with compositions and units. The {\em homology
$H^*(M)$} of a dg module $M$ is the induced functor
\[
H^*(\ca) \to \cg(k) \ko X \mapsto H^*(M(X))
\]
with values in the category $\cg(k)$ of graded $k$-modules
(\cf~\ref{ss:notations}). For each object $X$ of $\ca$, we have
the right module {\em represented by $X$}
\[
X^\wedge=\ca(?,X).
\]
The {\em category of dg modules $\cc(\ca)$} has as objects the dg
$\ca$-modules and as morphisms $L\to M$ the morphisms of dg
functors (\cf~\ref{ss:dgcat}). Note that $\cc(\ca)$ is an abelian
category and that a morphism $L\to M$ is an epimorphism
(respectively a monomorphism) iff it induces surjections
(respectively injections) in each component of $L(X) \to M(X)$ for
each object $X$ of $\ca$. A morphism $f: L \to M$ is a {\em
quasi-isomorphism} if it induces an isomorphism in homology.

We have $\cc(\ca)=Z^0(\cc_{dg}(\ca))$, where, in the notations
of~\ref{ss:dgcat}, the dg category $\cc_{dg}(\ca)$ is defined by
\[
\cc_{dg}(\ca)=\HOM(\ca\op, \cc_{dg}(k)).
\]
We write $\HOM(L,M)$ for the complex of morphisms from $L$ to $M$
in $\cc_{dg}(\ca)$. For each $X\in\ca$, we have a natural
isomorphism
\begin{equation} \label{eq:isorepgraded}
\HOM(X^\wedge,M) \iso M(X).
\end{equation}
The {\em category up to homotopy of dg $\ca$-modules} is
\[
\ch(\ca) = H^0(\cc_{dg}(\ca)).
\]
The isomorphism~\myref{eq:isorepgraded} yields isomorphisms
\begin{equation} \label{eq:isorephtpy}
\ch(\ca)(X^\wedge, M[n]) \iso H^n(\HOM(X^\wedge, M)) \iso H^n M(X)
\ko
\end{equation}
where $n\in\Z$ and $M[n]$ is the {\em shifted dg module} $Y
\mapsto M(Y)[n]$.

If $\ca$ is the dg category with one object associated with a
$k$-algebra $B$, then a dg $\ca$-module is the same as a complex
of $B$-modules.  More precisely, we have $\cc(\ca)=\cc(B)$,
$\cc_{dg}(\ca)=\cc_{dg}(B)$ and $\ch(\ca)=\ch(B)$.  In this case,
if $X$ is the unique object of $\ca$, the dg module $X^\wedge$ is
the complex formed by the free right $B$-module of rank one
concentrated in degree $0$.

\subsection{The derived category, resolutions}
The {\em derived category $\cd(\ca)$} is the localization of the
category $\cc(\ca)$ with respect to the class of
quasi-isomorphisms. Thus, its objects are the dg modules and its
morphisms are obtained from morphisms of dg modules by formally
inverting \cite{GabrielZisman67} all quasi-isomorphisms. The
projection functor $\cc(\ca) \to \cd(\ca)$ induces a functor
$\ch(\ca) \to \cd(\ca)$ and the derived category could
equivalently be defined as the localization of $\ch(\ca)$ with
respect to the class of all quasi-isomorphisms. Note that from
this definition, it is not clear that the morphism classes of
$\cd(\ca)$ are sets or that $\cd(\ca)$ is an additive category.

Call a dg module $P$ {\em cofibrant} if, for every surjective
quasi-isomorphism $L \to M$, every morphism $P \to M$ factors
through $L$. For example, for an object $X$  of $\ca$, the dg
module $X^\wedge$ is cofibrant. Call a dg module $I$ {\em fibrant}
if, for every injective quasi-isomorphism $L \to M$, every
morphism $L \to I$ extends to $M$. For example, if $E$ is an
injective cogenerator of the category of $k$-modules and $X$ an
object of $\ca$, the dg module $\HOM(\ca(X,?),E)$ is fibrant.

\begin{proposition}
\begin{itemize}
\item[a)] For each dg module $M$, there is a quasi-isomorphism
$\bp M \to M$ with cofibrant $\bp M$ and a quasi-isomorphism $M
\to \bi M$ with fibrant $\bi M$.
\item[b)] The projection functor $\ch(\ca) \to \cd(\ca)$ admits
a fully faithful left adjoint given by $M \mapsto \bp M$ and a
fully faithful right adjoint given by $M \mapsto \bi M$.
\end{itemize}
\end{proposition}

One can construct $\bp M$ and $\bi M$ explicitly, as first done in
\cite{AvramovHalperin86} (\cf also \cite{Keller94}). We call $\bp M
\to M$ a {\em cofibrant resolution} and $M \to \bi M$ a {\em
fibrant resolution} of $M$. According to b), these resolutions are
functorial in the category up to homotopy $\ch(\ca)$ and we can
compute morphisms in $\cd(\ca)$ via
\[
\ch(\ca)(\bp L, M) = \cd(\ca)(L,M) = \ch(\ca)(L, \bi M).
\]
In particular, for an object $X$ of $\ca$ and a dg module $M$, the
isomorphisms~\myref{eq:isorephtpy} yield
\begin{equation} \label{eq:isorepder}
\cd(\ca)(X^\wedge, M[n]) \iso \ch(\ca)(X^\wedge, M[n]) \iso H^n
M(X)
\end{equation}
since $X^\wedge$ is cofibrant. The embedding $\cd(\ca)\to
\ch(\ca)$ provided by $\bp$ also shows that the derived category
is additive.

If $\ca$ is associated with a $k$-algebra $B$ and $M$ is a right
$B$-module considered as a complex concentrated in degree $0$,
then $\bp M \to M$ is a projective resolution of $M$ and $M \to
\bi M$ an injective resolution. The proposition is best understood
in the language of Quillen model categories \cite{Quillen67}. We
refer to \cite{DwyerSpalinski95} for a highly readable
introduction and to \cite{Hovey99} \cite{Hirschhorn03} for
in-depth treatments. The proposition results from the following
theorem, proved using the techniques of \cite[2.3]{Hovey99}.

\begin{theorem} The category $\cc(\ca)$ admits two structures
of Quillen model category whose weak equivalences are the
quasi-isomorphisms:
\begin{itemize}
\item[1)] The {\em projective} structure, whose fibrations
are the epimorphisms. For this structure, each object is fibrant
and an object is cofibrant iff it is a cofibrant dg module.
\item[2)] The {\em injective} structure, whose cofibrations
are the monomorphisms. For this structure, each object is
cofibrant and an object is fibrant iff it is a fibrant dg module.
\end{itemize}
For both structures, two morphisms are homotopic iff they become
equal in the category up to homotopy $\ch(\ca)$.
\end{theorem}

\subsection{Exact categories, Frobenius categories}
Recall that an {\em exact category} in the sense of Quillen
\cite{Quillen73} is an additive category $\ce$ endowed with a
distinguished class of sequences
\[
\xymatrix{ 0 \ar[r] & A \ar[r]^i & B \ar[r]^p & C \ar[r] & 0} \ko
\]
where $i$ is a kernel of $p$ and $p$ a cokernel of $i$. We will
state the axioms these sequences have to satisfy using the
terminology of \cite{GabrielRoiter92}: The morphisms $p$ are
called deflations, the morphisms $i$ inflations and the pairs
$(i,p)$ conflations. The axioms are:
\begin{itemize}
\item[Ex0] The identity morphism of the zero object is a
deflation.
\item[Ex1] The composition of two deflations is a deflation.
\item[Ex2] Deflations are stable under base change.
\item[Ex2'] Inflations are stable under cobase change.
\end{itemize}
As shown in \cite{Keller90}, these axioms are equivalent to
Quillen's and they imply that if $\ce$ is small, then there is a
fully faithful functor from $\ce$ into an abelian category $\ce'$
whose image is an additive subcategory closed under extensions and
such that a sequence of $\ce$ is a conflation iff its image is a
short exact sequence of $\ce'$.  Conversely, one easily checks
that an extension closed full additive subcategory $\ce$ of an
abelian category $\ce'$ endowed with all conflations which become
exact sequences in $\ce'$ is always exact.  The fundamental
notions and constructions of homological algebra, and in
particular the construction of the derived category, naturally
extend from abelian to exact categories, \cf \cite{Neeman90} and
\cite{Keller96}.

A {\em Frobenius category} is an exact category $\ce$ which has
enough injectives and enough projectives and where the class of
projectives coincides with the class of injectives. In this case,
the {\em stable category} $\ul{\ce}$ obtained by dividing $\ce$ by
the ideal of morphisms factoring through a projective-injective
carries a canonical structure of triangulated category, \cf
\cite{Heller68} \cite{Happel87} \cite{KellerVossieck87}
\cite{GelfandManin96}. We write $\ol{f}$ for the image in
$\ul{\ce}$ of a morphism $f$ of $\ce$. The suspension functor $S$
of $\ul{\ce}$ is obtained by choosing a conflation
\[
\xymatrix{ 0 \ar[r] & A \ar[r] & IA \ar[r] & SA \ar[r] & 0}
\]
for each object $A$. Each triangle is isomorphic to a standard
triangle $(\ol{\imath}, \ol{p}, \ol{e})$ obtained by embedding a
conflation $(i,p)$ into a commutative diagram
\[
\xymatrix{ 0 \ar[r] & A \ar[r]^i \ar[d]_\id & B \ar[d] \ar[r]^{p}
&
                                          C \ar[d]^e \ar[r] & 0 \\
0 \ar[r] & A \ar[r] & IA \ar[r] & SA \ar[r] & 0.}
\]

\subsection{Triangulated structure} \label{ss:triangulated-structure}
Let $\ca$ be a small dg category. Define a sequence
\[
\xymatrix{0 \ar[r] & L \ar[r]^i & M \ar[r]^p & N \ar[r] & 0}
\]
of $\cc(\ca)$ to be a {\em conflation} if there is a morphism $r
\in \HOM(M,L)^0$ such that $ri=\id_L$ or, equivalently, a morphism
$s\in\HOM(N,M)$ such that $ps=\id_N$.

\begin{lemma}
\begin{itemize}
\item[a)] Endowed with these conflations, $\cc(\ca)$ becomes
a Frobenius category. The resulting stable category is canonically
isomorphic to $\ch(\ca)$. The suspension functor is induced by the
shift $M \mapsto M[1]$.
\item[b)] Endowed with the suspension induced by that of $\ch(\ca)$
and the triangles isomorphic to images of triangles of $\ch(\ca)$
the derived category $\cd(\ca)$ becomes a triangulated category.
Each short exact sequence of complexes yields a canonical
triangle.
\end{itemize}
\end{lemma}

\subsection{Compact objects, Brown representability}
\label{ss:CompactObjectsBrownRep} Let $\ct$ be a triangulated
category admitting arbitrary coproducts. Since the adjoint of a
triangle functor is a triangle functor \cite{KellerVossieck87},
the coproduct of triangles is then automatically a triangle.
Moreover, $\ct$ is {\em idempotent complete}
\cite{BoekstedtNeeman93}, \ie each idempotent endomorphism of an
object of $\ct$ is the composition of a section with a retraction.
An object $C$ of $\ct$ is {\em compact} if the functor $\ct(C,?)$
commutes with arbitrary coproducts, \ie for each family $(X_i)$ of
objects of $\ct$, the canonical morphism
\[
\coprod \ct(C,X_i) \to \ct(C, \coprod X_i)
\]
is invertible. The triangulated category $\ct$ is {\em compactly
generated} if there is a set $\cg$ of compact objects $G$ such
that an object $X$ of $\ct$ vanishes iff we have $\ct(G,X)=0$ for
each $G\in\cg$.

\begin{theorem}[Characterization of compact objects
\cite{ThomasonTrobaugh90} \cite{Neeman92a}] An object of $\ct$ is
compact iff it is a direct factor of an iterated extension of
copies of objects of $\cg$ shifted in both directions.
\end{theorem}

\begin{theorem}[Brown representability
\cite{Brown62} \cite{Adams71} \cite{Neeman96}] If $\ct$ is
compactly generated, a cohomological functor $F: \ct\op\to \Mod\Z$
is representable iff it takes coproducts of $\ct$ to products of
$\Mod\Z$.
\end{theorem}

A set of objects $\cg$ {\em symmetrically generates $\ct$} \cite{Krause02}
if we have
\begin{itemize}
\item[1)] an object $X$ of $\ct$ vanishes iff $\ct(G,X)=0$ for
each $G\in\cg$ and
\item[2)] there is a set of objects $\cg'$ such that
a morphism $f: X\to Y$ of $\ct$ induces surjections
$\ct(G,X)\to\ct(G,Y)$ for all $G\in\cg$ iff it induces injections
$\ct(Y,G') \to \ct(X,G')$ for all $G'\in\cg'$.
\end{itemize}
If $\cg$ compactly generates $\ct$, then we can take for $\cg'$
the set of objects $G'$ defined by
\[
\ct(?,G') = \Hom_k(\ct(G,?),E) \ko G\in\cg \ko
\]
where $E$ is an injective cogenerator of the category of
$k$-modules. Thus, in this case, $\cg$ also symmetrically generates $\ct$.

\begin{theorem}[Brown representability for the dual
\cite{Neeman98b} \cite{Krause02}] If $\ct$ is symmetrically
generated, a homological functor $F:\ct\to\Mod\Z$ is
corepresentable iff it commutes with products.
\end{theorem}

Let $\ca$ be a small dg category. The derived category $\cd(\ca)$
admits arbitrary coproducts and these are induced by coproducts of
modules. Thanks to the isomorphisms
\begin{equation} \label{eq:isorepgen}
\cd(\ca)(X^\wedge[n], M) \iso H^{-n} M(X)
\end{equation}
obtained from~\myref{eq:isorepder}, each dg module $X^\wedge[n]$,
where $X$ is an object of $\ca$ and $n$ an integer, is compact.
The isomorphism~\myref{eq:isorepgen} also shows that a dg module
$M$ vanishes in $\cd(\ca)$ iff each morphism $X^\wedge[n] \to M$
vanishes. Thus the set $\cg$ formed by the $X^\wedge[n]$,
$X\in\ca$, $n\in\Z$, is a set of compact generators for
$\cd(\ca)$. The {\em triangulated category $\tria(\ca)$ associated
with $\ca$} is the closure in $\cd(\ca)$ of the set of
representable functors $X^\wedge$, $X\in\ca$, under shifts in both
directions and extensions. The {\em category of perfect objects
$\per(\ca)$} the closure of $\tria(\ca)$ under passage to direct
factors in $\cd(\ca)$. The above theorems yield the

\begin{corollary} An object of $\cd(\ca)$ is compact iff it
lies in $\per(\ca)$. A cohomological functor $\cd(\ca)\op\to\Mod
k$ is representable iff it takes coproducts of $\cd(\ca)$ to
products of $\Mod k$. A homological functor $\cd(\ca)\to\Mod k$ is
corepresentable iff it commutes with products.
\end{corollary}

\subsection{Algebraic triangulated categories}
\label{ss:algtriacat} Let $\ct$ be a $k$-linear triangulated
category. We say that $\ct$ is {\em algebraic} if it is triangle
equivalent to $\ul{\ce}$ for some $k$-linear Frobenius category
$\ce$. It is easy to see that each $k$-linear triangulated
subcategory of an algebraic triangulated category is algebraic. We
will see below that each Verdier localization of an algebraic
triangulated category is algebraic (if we neglect a set-theoretic
problem). Moreover, categories of complexes up to homotopy are
algebraic, by~\ref{ss:triangulated-structure}. Therefore, `all'
triangulated categories occurring in algebra and geometry are
algebraic. Non algebraic  triangulated categories appear naturally
in topology (\cf also section~\ref{ss:topological}): For instance,
in the homotopy category of $2$-local spectra, the identity
morphism of the cone over twice the identity of the sphere
spectrum is of order four, but in each algebraic triangulated
category, the identity of the cone on twice the identity of an
object is of order two at most. A general method to prove that a
triangulated category obtained from a suitable stable Quillen
model category is not algebraic is to show that its
\cite{RezkSchwedeShipley01} Hom-functor enriched in spectra does
not factor through the canonical functor from the derived category
of abelian groups to the homotopy category of spectra.

We wish to show that `all' algebraic triangulated categories can
be described by dg categories. Let $\ct$ be a triangulated
category and $\cg$ a full subcategory. We make $\cg$ into a graded
category $\cg_{gr}$ by defining
\[
\cg_{gr}(G,G')= \bigoplus_{n\in\Z} \ct(G, S^n G').
\]
We obtain a natural functor $\ol{F}$ from $\ct$ to the category of
graded $\cg_{gr}$-modules by sending an object $Y$ of $\ct$ to
the $\cg_{gr}$-module
\[
X \mapsto \bigoplus_{n\in\Z}\ct(X,S^n Y)
\]

\begin{theorem}[\cite{Keller94}] Suppose that $\ct$ is algebraic. Then
there is a dg category $\ca$ such that $H^*(\ca)$ is isomorphic to
$\cg_{gr}$ and a triangle functor
\[
F: \ct \to \cd(\ca)
\]
such that the composition $H^*\circ F$ is isomorphic to $\ol{F}$.
Moreover,
\begin{itemize}
\item[a)] $F$ induces an equivalence
from $\ct$ to $\tria(\ca)$ iff $\ct$ coincides with its smallest
full triangulated subcategory containing $\cg$;
\item[b)] $F$ induces an equivalence from $\ct$ to $\per(\ca)$
iff $\ct$ is idempotent complete (\cf
section~\ref{ss:CompactObjectsBrownRep}) and equals the closure of
$\cg$ under shifts in both directions, extensions and passage to
direct factors;
\item[c)] $F$ is an equivalence $\ct \iso \cd(\ca)$
iff $\ct$ admits arbitrary coproducts and the objects of $\cg$
form a set of compact generators for $\ct$.
\end{itemize}
\end{theorem}

Examples arise from commutative and non commutative geometry:
A.~Bondal and M.~Van den Bergh show in \cite{BondalVandenBergh03}
that if $X$ is a quasi-compact quasi-separated scheme, then the
(unbounded) derived category $\ct=\cd_{qc}(X)$ of complexes of
$\co_X$-modules with quasi-coherent homology admits a single
compact generator $G$ and that moreover, $\Hom(G,G[n])$ vanishes
except for finitely many $n$. Thus $\ct$ is equivalent to the
derived category of a dg category with one object whose
endomorphism ring has bounded homology.

R.~Rouquier shows in \cite{Rouquier03} (\cf also
\cite{KrauseKussin05}) that if $X$ is a quasi-projective scheme
over a perfect field $k$, then the derived category of coherent
sheaves over $X$ admits a generator as a triangulated category (as
in part b) and, surprisingly, that it is even of `finite
dimension' as a triangulated category: each object occurs as a
direct factor of an object which admits a `resolution' of bounded
length by finite sums of shifts of the generator. Thus, the
bounded derived category of coherent sheaves is equivalent to
$\per(\ca)$ for a dg category with one object whose endomorphism
ring satisfies a strong regularity condition.

In \cite{Block05}, J.~Block describes the bounded derived category
of complexes of sheaves with coherent homology on a complex
manifold $X$ as the category $H^0(\ca)$ associated with a dg
category constructed from the Dolbeault dg algebra
$(A^{0,\bullet}(X), \ol{\partial})$. This can be seen as an
instance of a), where, for $\cg$, we can take for example the
category of coherent sheaves (\ie complexes concentrated in degree
0). Note however that the term `perfect derived category' is used
with a different meaning in \cite{Block05}.

In the independently obtained \cite{DwyerGreenlees02}, W.~Dwyer and
J.~Greenlees give elegant descriptions via dg endomorphism rings of
categories of complete, respectively torsion, modules. Their results
are applied in a unifying study of duality phenomena in algebra and
topology in \cite{DwyerGreenleesIyengarPrep05}.

One of the original motivations for the theorem was D.~Happel's
description \cite{Happel87} \cite{Happel88} of the bounded derived
category of a finite-dimensional associative algebra of finite
global dimension as the stable category of a certain Frobenius
category. This in turn was inspired by Bernstein-Gelfand-Gelfand's
\cite{BernsteinGelfandGelfand78} and Beilinson's
\cite{Beilinson78} descriptions of the derived category of
coherent sheaves on projective space.

A vast generalization of the theorem to non-additive contexts
\cite{SchwedeShipley00} is due to S.~Schwede and B.~Shipley
\cite{SchwedeShipley03}, \cf also section~\ref{ss:topological}
below.

\subsection{Well-generated algebraic triangulated categories}
A triangulated category $\ct$ is {\em well-generated}
\cite{Neeman99} \cite{Krause01} if it admits arbitrary coproducts
and a {\em good set of generators $\cg$}, \ie $\cg$ is stable
under shifts in both directions and satisfies
\begin{itemize}
\item[1)] an object $X$ of $\ct$ vanishes iff
$\ct(G,X)=0$ for each $G\in\cg$,
\item[2)] there is a cardinal $\alpha$ such that each $G\in\cg$ is
{\em $\alpha$-small}, \ie for each family of objects $X_i$,
$i\in I$, of $\ct$, each morphism
\[
G \to \bigoplus_{i\in I} X_i
\]
factors through a subsum $\bigoplus_{i\in J} X_i$ for some subset
$J$ of $I$ of cardinality strictly smaller than $\alpha$,
\item[3)] for each family of morphisms $f_i: X_i \to Y_i$, $i\in I$,
of $\ct$ which induces surjections
\[
\ct(G,X_i) \to \ct(G,Y_i)
\]
for all $G\in \cg$ and all $i\in I$, the sum of the $f_i$ induces
surjections
\[
\ct(G,\bigoplus X_i) \to \ct(G,\bigoplus Y_i)
\]
for all $G\in\cg$.
\end{itemize}
Clearly each compactly generated triangulated category is
well-generated. A.~Neeman proves in \cite{Neeman99} that the Brown
representability theorem holds for well-genera\-ted triangulated
categories. This is one of the main reasons for studying them.
Another important result of \cite{Neeman99} is that if $\ct$ is
well-generated and $\cs\to \ct$ is a localization (\ie a fully
faithful triangle functor admitting a left adjoint whose kernel
is generated by a {\em set} of objects) then $\cs$ is
well-generated. Thus each localization of a compactly generated
triangulated category is well-generated and in particular, so is
each localization of the derived category of a small dg category.

Here is another class of examples: Let $\cb$ be a Grothendieck
abelian category, \eg the category of modules on a ringed space.
Then, by the Popescu-Gabriel theorem \cite{PopescuGabriel64}
\cite{Lowen04a}, $\cb$ is the localization of the category of
$\Mod A$ of $A$-modules over some ring $A$. One can
deduce from this that the unbounded derived category of
of the abelian category $\cb$ (\cf \cite{Franke97}
\cite{Alonso-Jeremias-Souto00} \cite{KashiwaraSchapira05}) is a
localization of $\cd(A)$ and thus is well-generated.

\begin{theorem}[\cite{Porta06}] Let $\ct$ be an algebraic triangulated
category. Then $\ct$ is well-generated iff it is a localization of
$\cd(\ca)$ for some small dg category $\ca$. Moreover, if $\ct$ is
well-generated and $\cu\subset\ct$ a full small subcategory such
that, for each $X\in\ct$, we have
\[
X=0 \Leftrightarrow \ct(U,S^n X)=0 \mbox{ for all } n\in\Z \mbox{
and } U\in \cu \ko
\]
then there is an associated localization $ \ct \to \cd(\ca) $ for
some small dg category $\ca$ with $H^*(\ca)=\cu_{gr}$.
\end{theorem}

\subsection{Morita equivalence} \label{ss:Morita-equivalence}
Let $\ca$ and $\cb$ be small dg categories. Let $X$ be an {\em
$\ca$-$\cb$-bimodule}, \ie a dg $\ca\op\ten\cb$-module $X$. Thus
$X$ is given by complexes $X(B,A)$, for all $A$ in $\ca$ and $B$
in $\cb$, and morphisms of complexes
\[
\cb(A,A')\ten X(B,A) \ten \ca(B',B) \to X(B',A').
\]
For each dg $\cb$-module $M$, we obtain a dg $\ca$-module
\[
GM=\HOM(X,M) : A \mapsto \HOM(X(?,A), M).
\]
The functor $G: \cc(\cb) \to \cc(\ca)$ admits a left adjoint $F:L
\mapsto L\ten_\ca X$. These functors do not respect
quasi-isomorphisms in general, but they form a Quillen adjunction
(\cf section~\ref{ss:topological}) and their derived functors
\[
\L F : L \mapsto F(\bp L) \mbox{ and } \R G : M \mapsto G(\bi M)
\]
form an adjoint pair of functors between $\cd(\ca)$ and
$\cd(\cb)$.

\begin{lemma}[\cite{Keller94}] \label{lemma:tensor-equivalence}
The functor $\L F: \cd(\ca) \to \cd(\cb)$ is an equivalence if and
only if
\begin{itemize}
\item[a)] the dg $\cb$-module $X(?,A)$ is perfect for all $A$ in $\ca$,
\item[b)] the morphism
\[
\ca(A,A') \to \HOM(X(?,A), X(?,A'))
\]
is a quasi-isomorphism for all $A$, $A'$ in $\ca$ and
\item[c)] the dg $\cb$-modules $X(?,A)$, $A\in\ca$, form a set
of (compact) generators for $\cd(\cb)$.
\end{itemize}
\end{lemma}

For example, if $E: \ca \to \cb$ is a dg functor, then
$X(B,A)=\cb(B,E(A))$ defines a dg bimodule so that the above
functor $G$ is the restriction along $E$. Then the lemma shows
that $\R G$ is an equivalence iff $E$ is a quasi-equivalence. We
loosely refer to the functor $\L F$ associated with a dg
$\ca$-$\cb$-bimodule as a {\em tensor functor}.

\begin{theorem}[\cite{Keller94}] The following are equivalent
\item[1)] There is an equivalence $\cd(\ca) \to \cd(\cb)$ given
by a composition of tensor functors and their inverses.
\item[2)] There is a dg subcategory $\cg$ of $\cc(\cb)$ formed by
cofibrant dg modules such that the objects of $\cg$ form a set of
compact generators for $\cd(\cb)$ and there is a chain of
quasi-equivalences
\[
\ca \la \ca' \ra \cdots \la \cg' \ra \cg
\]
linking $\ca$ to $\cg$.
\end{theorem}

We say that $\ca$ and $\cb$ are {\em dg Morita equivalent} if the
conditions of the theorem are satisfied. In this case, there is of
course a triangle equivalence $\cd(\ca)\to \cd(\cb)$. In general,
the existence of such a triangle equivalence is not sufficient for
$\ca$ and $\cb$ to be dg Morita equivalent, \cf
section~\ref{ss:topological}. The following theorem is therefore
remarkable:

\begin{theorem}[Rickard \cite{Rickard89}] \label{thm:Rickard}
Suppose that $\ca$ and $\cb$ have their homology concentrated in
degree $0$. Then the following are equivalent:
\begin{itemize}
\item[1)] $\ca$ and $\cb$ are dg Morita equivalent.
\item[2)] There is a triangle equivalence $\cd(\ca)\to\cd(\cb)$.
\item[3)] There is a full subcategory $\ct$ of $\cd(\cb)$ such that
\begin{itemize}
\item[a)] the objects of $\ct$ form a set of compact generators of $\cd(\cb)$,
\item[b)] we have $\cd(\cb)(T,T'[n])=0$ for all $n\neq 0$ and all $T$, $T'$
of $\ct$,
\item[c)] there is an equivalence $H^0(\ca) \iso \ct$.
\end{itemize}
\end{itemize}
\end{theorem}

We refer to \cite{Keller94} or \cite{DuggerShipley04} for this
form of the theorem. A subcategory $\ct$ satisfying a) and b) in
3) is called a {\em tilting subcategory}, a concept which
generalizes that of a tilting module. We refer to \cite{Reiten98}
\cite{AngeleriHappelKrause06} for the theory of tilting, from
which this theorem arose and which provides huge classes of
examples from the representation theory of finite-dimensional
algebras and finite groups as well as from algebraic geometry, \cf
also the appendix to \cite{Neeman05} and \cite{KoenigZimmermann98}
\cite{Rickard98} \cite{Rouquier06}.

\subsection{Topological Morita equivalence} \label{ss:topological}
In recent years, Morita theory has been vastly generalized from
algebraic triangulated categories to stable model categories in
work due to S.~Schwede and B.~Shipley. This is based on the
category of symmetric spectra as constructed in
\cite{HoveyShipleySmith00} (\cf \cite{EKMM97} for a different
construction of a symmetric monoidal model category for the
category of spectra).  We refer to \cite{Schwede04} for an
excellent exposition of these far-reaching results and their
surprising applications in homotopy theory.  In work by D.~Dugger
and B.~Shipley \cite{DuggerShipley06a}, \cf also
\cite{DuggerShipley06b} \cite{DuggerShipley06c}, this `topological
Morita theory' has been applied to dg categories. We briefly
describe their results and refer to \cite{Shipley06} for a highly
readable, more detailed survey.

The main idea is to replace the monoidal base category, the
derived category of abelian groups $\cd(\Z)$, by a more
fundamental category: the `derived category of the category of
sets', \ie the homotopy category of spectra. To preserve higher
homotopical information, one must not, of course, work at the
level of derived categories but has to introduce model categories.
So instead of considering $\cd(\Z)$, one considers its model
category $\cc(\Z)$ of complexes of abelian groups and replaces it
by a convenient model of the category of spectra: the category of
symmetric spectra, which one might imagine as `complexes of
abelian groups up to homotopy'. We refer to
\cite{HoveyShipleySmith00} or \cite{Schwede04} for the precise
definition. As shown in \cite{HoveyShipleySmith00}, symmetric
spectra form a {\em symmetric} monoidal category which carries a
compatible Quillen model structure and whose homotopy category is
equivalent to the homotopy category of spectra of Bousfield and
Friedlander \cite{BousfieldFriedlander78}. The tensor product is
the {\em smash product $\wedge$} and the unit object is the {\em
sphere spectrum $\S$}. The unit object is cofibrant and the smash
product induces a monoidal structure on the homotopy category of
symmetric spectra.  The {\em Eilenberg-MacLane functor $H$} is a
lax monoidal functor from the category of complexes $\cc(k)$ to
the category of symmetric spectra such that the homology groups of
a complex $C$ become isomorphic to the homotopy groups of $HC$.
Since $H$ is lax monoidal, if $A$ is a dg $\Z$-algebra, then $HA$
is naturally an algebra in the category of symmetric spectra and
if $M$ is an $A$-module, then $HM$ becomes an $HA$-module.  More
generally, if $\ca$ is a dg category over $\Z$, then $H\ca$
becomes a {\em spectral category}, \ie a category enriched in
symmetric spectra, \cf \cite{SchwedeShipley03b} \cite{Shipley02}.
Each $\ca$-module $M$ then gives rise to a {\em spectral module
$HM$} over $H\ca$. The spectral modules over a spectral category
form a Quillen model category \cite{Shipley02}.

Recall that if $\cl$ and $\cm$ are Quillen model categories, a
{\em Quillen adjunction} is given by a pair of adjoint functors
$L:\cl\to\cm$ and $R:\cm\to\cl$ such that $L$ preserves
cofibrations and $R$ fibrations. Such a pair induces an adjoint
pair between the homotopy categories of $\cl$ and $\cm$. If the
induced functors are equivalences, then $(L,R)$ is a {\em Quillen
equivalence}. The model categories $\cl$ and $\cm$ are {\em
Quillen equivalent} if they are linked by a chain of Quillen
equivalences.

It was shown by A.~Robinson \cite{Robinson87}, \cf also
\cite{SchwedeShipley03}, that for an ordinary ring $R$, the
unbounded derived category of $R$-modules is equivalent to the
homotopy category of spectral modules over $HR$. This result is
generalized and refined as follows:

\begin{theorem}[Shipley \cite{Shipley02}] If $\ca$ is a dg category over $\Z$,
the model categories of dg $\ca$-modules and of spectral modules
over $H\ca$ are Quillen equivalent.
\end{theorem}

This allows us to define two small dg categories $\ca$ and $\cb$
to be {\em topologically Morita equivalent} if their categories of
spectral modules are Quillen equivalent.

\begin{proposition}[\cite{DuggerShipley06a}]
Let $\ca$ and $\cb$ be two dg rings. Then statement a) implies b)
and b) implies c):
\begin{itemize}
\item[a)] $\ca$ and $\cb$ are dg Morita equivalent.
\item[b)] $\ca$ and $\cb$ are topologically Morita equivalent.
\item[c)] $\cd(\ca)$ is triangle equivalent to $\cd(\cb)$.
\end{itemize}
\end{proposition}

It is remarkable that in general, these implications are strict.
Examples which show this were obtained in recent joint work by
D.~Dugger and B.~Shipley \cite{DuggerShipley06a}, \cf also
\cite{Shipley06}. To show that c) does not imply b), they invoke
Schlichting's example \cite{Schlichting02}: Let $p$ be an odd
prime. The module categories over $A'=\Z/p^2$ and
$B'=(\Z/p)[\eps]/\eps^2$ are Frobenius categories. Their stable
categories are triangle equivalent (both are equivalent to the
category of $\Z/p$-vector spaces with the identical suspension and
the split triangles) but the $K$-theories associated with the
stable module categories are not isomorphic. Since $K$-theory is
preserved under topological Morita equivalence (\cf
section~\ref{ss:K-theory} below), the dg algebras $A$ and $B$
associated (\cf section~\ref{ss:algtriacat}) with the canonical
generators (corresponding to the one-dimensional vector space over
$\Z/p$) of the stable categories of $A'$ and $B'$ cannot be
topologically Morita equivalent.

To show that b) does not imply a), Dugger and Shipley consider two
dg algebras $A$ and $B$ with homology isomorphic to $\Z/2\oplus
\Z/2[2]$. The isomorphism classes of such algebras in the homotopy
category of dg $\Z$-algebras are parametrized by the Hochschild
cohomology group $HH^4_\Z(\Z/2,\Z/2)$. Their isomorphism classes
in the homotopy category of $\S$-algebras are parametrized by the
topological Hochschild cohomology group $THH^4_\S(\Z/2,\Z/2)$ as
shown in \cite{Lazarev01}. The computation of the Hochschild
cohomology group $HH^4_\Z(\Z/2,\Z/2)$ is elementary and, thanks to
Fran\-jou-Lannes-Schwartz' work \cite{FranjouLannesSchwartz94},
the topological Hochschild cohomology algebra
$$
THH^*_\S(\Z/2,\Z/2)
$$ is known. Dugger-Shipley then conclude by exhibiting a non-trivial
element in the kernel of the canonical map
\[
\Phi: HH^4_\Z(\Z/2,\Z/2) \to THH^4_\S(\Z/2,\Z/2).
\]
The explicit description of the two algebras is given in
\cite{Shipley02} \cite{DuggerShipley04} \cite{DuggerShipley06a}.
The appearance of torsion in these examples is unavoidable: for dg
algebras over the rationals, statements a) and b) above are
equivalent \cite{DuggerShipley06a}.

\section{The homotopy category of small dg categories}
\label{s:HtpyCatSmalldgCat}

\subsection{Introduction} \label{ss:IntroHtpyCat}
Invariants like $K$-theory, Hochschild homology, cyclic homology
\ldots na\-turally extend from $k$-algebras to dg categories (\cf
section~\ref{s:invariants}). In analogy with the case of ordinary
$k$-algebras, these extended invariants are preserved under dg
Morita equivalence. However, unlike the module category over a
$k$-algebra, the derived category of a dg category, even with its
triangulated structure, does not contain enough information to
compute the invariant (\cf the examples in
section~\ref{ss:topological}).  Our aim in this section is to
present a category obtained from that of small dg categories by
`inverting the dg Morita equivalences'. It could be called the
`homotopy category of enhanced (idempotent complete) triangulated
categories' \cite{BondalKapranov89} or the `Morita homotopy
category of small dg categories' $\Hmo$, as in \cite{Tabuada05b}.
Invariants like $K$-theory and cyclic homology factor through the
Morita homotopy category.

The Morita homotopy category very much resembles the category of
small, idempotent complete, triangulated categories. In
particular, it admits `dg quotients' \cite{Drinfeld04}, which
correspond to Verdier localizations. Like these, they are
characterized by a universal property. The great advantages of the
Morita homotopy category over that of small triangulated
categories are that moreover, it admits {\em all (homotopy) limits
and colimits} (like any homotopy category of a Quillen model
category) and is {\em monoidal and closed}.

The Morita homotopy category $\Hmo$ is a full subcategory of the
localization $\Hqe$ of the category of small dg categories with
respect to the quasi-equivalences. The first step is therefore to
analyze the larger category $\Hqe$. Its morphism spaces are
revealed by To\"en's theorem~\ref{thm:mapping-spaces} below.

\subsection{Inverting quasi-equivalences}
\label{ss:InvertingQuasiEq} Let $k$ be a commutative ring and
$\dgcat_k$ the category of small dg $k$-categories as in
section~\ref{ss:diffgradedcat}. An analogue of the following
theorem for simplicial categories is proved in \cite{Bergner04}.

\begin{theorem}[\cite{Tabuada05a}] \label{thm:Tabuada1}
The category $\dgcat_k$ admits a structure of cofibrantly
generated model category whose weak equivalences are the
quasi-equivalences and whose fibrations are the dg functors
$F:\ca\to\cb$ which induce componentwise surjections
$\ca(X,Y)\to\cb(FX,FY)$ for all $X,Y$ in $\ca$ and such that, for
each isomorphism $v: F(X) \to Z$ of $H^0(\cb)$, there is an
isomorphism $u$ of $H^0(\ca)$ with $F(u)=v$.
\end{theorem}

This shows in particular that the {\em localization $\Hqe$} of
$\dgcat_k$ with respect to the quasi-equivalences has small
$\Hom$-sets and that we can compute morphisms from $\ca$ to $\cb$
in the localization as morphisms modulo homotopy from a cofibrant
replacement $\ca_{cof}$ of $\ca$ to $\cb$ (note that all small dg
categories are fibrant). In general, the cofibrant replacement
$\ca_{cof}$ is not easy to compute with but if $\ca(X,Y)$ is
cofibrant in $\cc(k)$ and the unit morphisms $k\to \ca(X,X)$ admit
retractions in $\cc(k)$ for all objects $X$, $Y$ of $\ca$, for
example if $k$ is a field, then for $\ca_{cof}$, we can take the
category with the same objects as $\ca$ and whose morphism spaces
are given by the `reduced cobar-bar construction', \cf \eg
\cite{Drinfeld04} \cite{Keller05c}. The homotopy relation is then
the one of \cite[3.3]{Keller99}.

However, the morphism sets in the localization are much better
described as follows: Consider two dg categories $\ca$ and $\cb$.
If necessary, we replace $\ca$ by a quasi-equivalent dg category
so as to achieve that $\ca$ is {\em $k$-flat}, \ie the functor
$\ca(X,Y)\ten ?$ preserves quasi-isomorphisms for all $X$, $Y$ of
$\ca$ (for example, we could take a cofibrant replacement of
$\ca$). Let $\rep(\ca,\cb)$ be the full subcategory of the derived
category $\cd(\ca\op\ten\cb)$ of $\ca$-$\cb$-bimodules formed by
the bimodules $X$ such that the tensor functor
\[
?\lten_\ca X: \cd(\ca) \to \cd(\cb)
\]
takes the representable $\ca$-modules to objects which are
isomorphic to representable $\cb$-modules. In other words, we
require that $X(?,A)$ is quasi-isomorphic to a representable
$\cb$-module for each object $A$ of $\ca$. We call such a bimodule
a {\em quasi-functor} since it yields a genuine functor
\[
H^0(\ca) \to H^0(\cb).
\]
We think of $\rep(\ca,\cb)$ as the `category of representations up
to homotopy of $\ca$ in $\cb$'.

\begin{theorem}[To\"en \cite{Toen04}] The morphisms from $\ca$ to $\cb$
in the localization of $\dgcat_k$ with respect to the
quasi-equivalences are in natural bijection with the isomorphism
classes of $\rep(\ca,\cb)$.
\end{theorem}

The theorem has been in limbo for some time, \cf
\cite[2.3]{Keller98} \cite{Keller99} \cite{Drinfeld04}. It is due
to B.~To\"en, as a corollary of a much more precise statement:
Recall from \cite[Ch.~5]{Hovey99} that each model category $\cm$
admits a mapping space bifunctor
\[
\Map : \Ho(\cm)\op\times \Ho(\cm) \to \Ho(\Sset)
\]
such that we have, for example, the natural isomorphisms
\[
\pi_0(\Map(X,Y)) = \Ho(\cm)(X,Y).
\]
The spaces $\Map$ may also be viewed as the morphism spaces in the
Dwyer-Kan localization \cite{DwyerKan80a} \cite{DwyerKan80b}
\cite{DwyerKan80c} of $\cm$ with respect to the class of weak
equivalences, \cf \cite{DwyerKan80c} \cite{Hirschhorn03}. Now let
$\cR(\ca,\cb)$ be the category with the same objects as
$\rep(\ca,\cb)$ and whose morphisms are the quasi-isomorphisms of
dg bimodules. Thus, the category $\cR(\ca,\cb)$ is a non full
subcategory of the category of dg bimodules $\cc(\ca\op\ten\cb)$.

\begin{theorem}[To\"en \cite{Toen04}] \label{thm:mapping-spaces}
There is a canonical weak equivalence of simplicial sets between
$\Map(\ca,\cb)$ and the nerve of the category $\cR(\ca,\cb)$.
\end{theorem}

The theorem allows one to compute the homotopy groups of the {\em
classifying space $|\dgcat|$ of dg categories}, which is defined
as the nerve of the category of quasi-equivalences between dg
categories. Of course, the connected components of this space are
in bijection with the isomorphism classes of $\Hqe$. Now let $\ca$
be a small dg category. Then the fundamental group of $|\dgcat|$
at $\ca$ is the group of automorphisms of $\ca$ in $\Hqe$ (\cf
\cite{Quillen73}). For example, if $\ca$ is the category of
bounded complexes of projective $B$-modules over an ordinary
$k$-algebra $B$, then this group is the derived Picard group of
$B$ as studied in \cite{RouquierZimmermann03} \cite{Keller04}
\cite{Yekutieli04}. For the higher homotopy groups, we have the

\begin{corollary}[\cite{Toen04}]
\begin{itemize}
\item[a)] The group $\pi_2(|\dgcat|, \ca)$ is the group of
invertible elements of the dg center of $\ca$ (=its zeroth
Hochschild cohomology group).
\item[b)] For $i\geq 2$, the group $\pi_i(|\dgcat|, \ca)$ is the
$(2-i)$-th Hochschild cohomology of $\ca$.
\end{itemize}
\end{corollary}

\subsection{Closed monoidal structure}
\label{ss:ClosedMonoidalStructure} As we have observed in
section~\ref{ss:diffgradedcat}, the category $\dgcat_k$ admits a
tensor product $\ten$ and an internal $\Hom$-functor $\HOM$. If
$\ca$ is cofibrant, then the functor $\ca\ten ?$ preserves weak
equivalences so that the localization $\Hqe$ inherits a tensor
product $\lten$. However, the tensor product of two cofibrant dg
categories is not cofibrant in general (in analogy with the fact
that the tensor product of two non commutative free algebras is
not non commutative free in general). By the adjunction formula
\[
\HOM(\ca,\HOM(\cb,\cc))=\HOM(\ca\ten\cb,\cc) \ko
\]
it follows that even if $\ca$ is cofibrant, the functor
$\HOM(\ca,?)$ cannot preserve weak equivalences in general and
thus will not induce an internal $\Hom$-functor in $\Hqe$.
Nevertheless, we have the

\begin{theorem}[\cite{Drinfeld04} \cite{Toen04}] \label{thm:closed-monoidal}
The monoidal category $(\Hqe,\lten)$ admits an internal
$\Hom$-functor $\RHOM$. For two dg categories $\ca$ and $\cb$ such
that $\ca$ is $k$-flat, the dg category $\RHOM(\ca,\cb)$ is
isomorphic in $\Hqe$ to the dg category $\rep_{dg}(\ca,\cb)$, \ie
the full subcategory of the dg category of $\ca$-$\cb$-bimodules
whose objects are those of $\rep(\ca,\cb)$ and which are cofibrant
as bimodules.
\end{theorem}

Thus we have equivalences (we suppose $\ca$ $k$-flat)
\[
H^0(\RHOM(\ca,\cb))=H^0(\rep_{dg}(\ca,\cb)) \iso \rep(\ca,\cb).
\]
In terms of the internal $\Hom$-functor $\HOM$ of $\dgcat_k$, we
have
\[
H^0(\RHOM(\ca,\cb))= H^0(\HOM(\ca,\cb))[\Sigma^{-1}] \ko
\]
where $\Sigma$ is the set of morphisms $\phi: F \to G$ such that
$\phi(A)$ is invertible in $H^0(\cb)$ for all objects $A$ of
$\ca$, \cf \cite{Keller98}.

Yet another description can be given in terms of
$A_\infty$-functors: Let $\ca$ be a dg category such that the
morphism spaces $\ca(A,A')$ are cofibrant in $\cc(k)$ and the unit
maps $k \to \ca(A,A)$ admit retractions in $\cc(k)$ for all
objects $A$, $A'$ of $\ca$. Then the dg category $\RHOM(\ca,\cb)$
is quasi-equivalent to the $A_\infty$-category of (strictly
unital) $A_\infty$-functors from $\ca$ to $\cb$, \cf
\cite{Kontsevich98} \cite{Lefevre03} \cite{Lyubashenko03}
\cite{Keller05c}. Since $\cb$ is a dg category, this
$A_\infty$-category is in fact a dg category.

An important point of classical Morita theory is that for two
rings $B$, $C$, there is an equivalence between the category of
$B$-$C$-bimodules and the category of coproduct preserving
functors from the category of $B$-modules to that of $C$-modules
(note that here and in what follows, we need to consider `large'
categories and should introduce universes to make our statements
rigorous \ldots). Similarly, if $\ca$ is a small $k$-flat dg
category, we consider the large dg category $\cd_{dg}(\ca)$: it is
the full dg subcategory of $\cc_{dg}(\ca)$ whose objects are all
the cofibrant dg modules. Thus we have an equivalence of
categories
\[
\cd(\ca) = H^0(\cd_{dg}(\ca)).
\]
This shows that if $\cb$ is another dg category, then each
quasi-functor $X$ in
\[
\rep(\cd_{dg}(\ca),\cd_{dg}(\cb))
\]
gives rise to a functor $\cd(\ca) \to \cd(\cb)$. We say that the
quasifunctor $X$ {\em preserves coproducts} if this functor
preserves coproducts.

\begin{theorem}[\cite{Toen04}]
There is a canonical isomorphism in $\Hqe$
\[
\cd_{dg}(\ca\op\ten\cb) \iso \RHOM_c(\cd_{dg}(\ca), \cd_{dg}(\cb))
\ko
\]
where $\RHOM_c$ denotes the full subcategory of $\RHOM$ formed by
the coproduct preserving quasifunctors.
\end{theorem}

If we apply this theorem to $\cb=\ca$ and compare the endomorphism
algebras of the identity functors on both sides, we see that the
Hochschild cohomology (\cf section~\ref{ss:Hochschild-cohomology}
below) of the small dg category $\ca$ coincides with the
Hochschild cohomology of the large dg category $\cd_{dg}(\ca)$,
which is quite surprising. An analogous result for Grothendieck
abelian categories (in particular, module categories) is due to
T.~Lowen and M.~Van den Bergh \cite{LowenVandenBergh04a}.

\subsection{Dg localizations, dg quotients, dg-derived categories}
\label{ss:DgLocalizations}
 Let $\ca$
be a small dg category. Let $S$ be a set of morphisms of
$H^0(\ca)$. Let us say that a morphism $R:\ca\to\cb$ of $\Hqe$
{\em makes $S$ invertible} if the induced functor
\[
H^0(\ca) \to H^0(\cb)
\]
takes each $s\in S$ to an isomorphism.

\begin{theorem}[\cite{Toen04}] There is a morphism $Q: \ca\to \ca[S^{-1}]$ of
$\Hqe$ such that $Q$ makes $S$ invertible and each morphism $R$ of
$\Hqe$ which makes $S$ invertible uniquely factors through $Q$.
\end{theorem}

We call $\ca[S^{-1}]$ the {\em dg localization of $\ca$ at $S$}.
Note that it is unique up to unique isomorphism in $\Hqe$. It is
constructed in \cite{Toen04} as a homotopy pushout
\[
\xymatrix{ \coprod_{s\in S} I \ar[r] \ar[d] & \ca \ar[d] \\
\coprod_{s\in S} k \ar[r] & \ca[S^{-1}], }
\]
where $I$ denotes the dg $k$-category freely generated by one
arrow $f: 0 \to 1$ of degree $0$ with $df=0$ and left vertical
arrow is induced by the morphisms $I\to k$ which sends $f$ to $1$.
The universal property of $Q: \ca\to\ca[S^{-1}]$ admits refined
forms, namely, $Q$ induces an equivalence of categories
\[
\rep(\ca[S^{-1}], \cb) \iso \rep_S(\ca,\cb) \ko
\]
an isomorphism of $\Hqe$
\[
\rep_{dg}(\ca[S^{-1}],\cb) \iso \rep_{dg,S}(\ca,\cb) \ko
\]
and a weak equivalence of simplicial sets
\[
\Map(\ca[S^{-1}], \cb) \iso \Map_S(\ca,\cb).
\]
Here $\rep_S$ and $\rep_{dg,S}$ denote the full subcategories of
quasi-functors whose associated functors $H^0(\ca)\to H^0(\cb)$
make $S$ invertible and $\Map_S$ the union of the connected
components containing these quasi-functors.

An important variant is the following: Let $\cn$ be a set of
objects of $\ca$. Let us say that a morphism $Q :\ca\to\cb$ of
$\Hqe$ {\em annihilates $\cn$} if the induced functor
\[
H^0(\ca) \to H^0(\cb)
\]
takes all objects of $\cn$ to zero objects (\ie objects whose
identity morphism vanishes in $H^0(\cb)$).

\begin{theorem}[\cite{Keller99} \cite{Drinfeld04}]
There is a morphism $Q: \ca\to\ca/\cn$ of $\Hqe$ which annihilates
$\cn$ and is universal among the morphisms annihilating $\cn$.
\end{theorem}

We call $\ca/\cn$ the {\em dg quotient of $\ca$ by $\cn$}. If
$\ca$ is $k$-flat (\cf section~\ref{ss:InvertingQuasiEq}), then
$\ca/\cn$ admits a  beautiful simple construction
\cite{Drinfeld04}: One adjoins to $\ca$ a contracting homotopy for
each object of $\cn$. The general case can be reduced to this one
or treated using orthogonal subcategories \cite{Keller99}. The dg
quotient has refined universal properties analogous to those of
the dg localization. In particular, the morphism $\ca\to\ca/\cn$
induces an equivalence \cite{Drinfeld04}
\[
\rep(\ca/\cn,\cb) \to \rep_\cn(\ca,\cb) \ko
\]
where $\rep_\cn$ denotes the full subcategory of quasi-functors
whose associated functors $H^0(\ca)\to H^0(\cb)$ annihilate $\cn$.

Dg quotients yield functorial dg versions of Verdier localizations
\cite{Verdier96}. For example, if $\ce$ is a small abelian (or,
more generally, exact) category, we can take for $\ca$ the dg
category of bounded complexes $\cc^b_{dg}(\ce)$ over $\ce$ and for
$\cn$ the dg category of bounded acyclic complexes
$\ac^b_{dg}(\ce)$. Then we obtain the {\em dg-derived category}
\[
\cd^b_{dg}(\ce) = \cc^b_{dg}(\ce)/\ac^b_{dg}(\ce)
\]
so that we have
\[
\cd^b(\ce) = H^0(\cd^b_{dg}(\ce)).
\]
More generally, every localization pair \cite{Keller99}
(=Frobenius pair \cite{Schlichting06}) gives rise to a dg
category. After taking the necessary set-theoretic precautions, we
also obtain a dg-derived category
\[
\cd_{dg}(\ce) = \cc_{dg}(\ce)/\ac_{dg}(\ce)
\]
which refines the {\em unbounded} derived category of a $k$-linear
Grothendieck abelian category $\ce$. For a quasi-compact
quasi-separated scheme $X$, let us write $\cd_{dg}(X)$ for
$\cd_{dg}(\ce)$, where $\ce$ is the Grothendieck abelian category
of quasi-coherent sheaves on $X$. The following theorem shows that
dg functors between dg derived categories are much more closely related to
geometry than triangle functors between derived categories,
\cf \cite{BondalOrlov02} \cite{Orlov97}.

\begin{theorem}[\cite{Toen04}] Let $X$ and $Y$ be quasi-compact
separated schemes over $k$ such that $X$ is flat over $\Spec k$.
Then we have a canonical isomorphism in $\Hqe$
\[
\cd_{dg}(X\times_k Y) \iso \RHOM_c(\cd_{dg}(X),\cd_{dg}(Y)) \ko
\]
where $\RHOM_c$ denotes the full subcategory of $\RHOM$ formed by
the coproduct preserving quasi-functors. Moreover, if $X$ and $Y$
are smooth and projective over $\Spec k$, we have a canonical
isomorphism in $\Hqe$
\[
\parf_{dg}(X\times_k Y) \iso \RHOM(\parf_{dg}(X),\parf_{dg}(Y))
\]
where $\parf_{dg}$ denotes the full dg subcategory of $\cd_{dg}$
whose objects are the perfect complexes.
\end{theorem}

\subsection{Pretriangulated dg categories}
\label{ss:Pretriangulated} Let $\ca$ be a small dg category. We
say that $\ca$ is {\em pretriangulated} or {\em exact} if the
image of the Yoneda functor
\[
Z^0(\ca) \to \cc(\ca) \ko X \mapsto X^\wedge
\]
is stable under shifts in both directions and extensions (in the
sense of the exact structure of
section~\ref{ss:triangulated-structure}). Equivalently, for all
objects $X$, $Y$ of $\ca$ and all integers $n$, the object
$X^\wedge[n]$ is isomorphic to $X[n]^\wedge$ and the cone over a
morphism $f^\wedge: X^\wedge \to Y^\wedge$ is isomorphic to
$C(f)^\wedge$ for unique objects $X[n]$ and $C(f)$ of $Z^0(\ca)$.
If $\ca$ is exact, then $Z^0(\ca)$ becomes a Frobenius subcategory
of $\cc(\ca)$ and $H^0(\ca)$ a triangulated subcategory of
$\ch(\ca)$. If $\cb$ is an exact dg category and $\ca$ an
arbitrary dg category, then $\HOM(\ca,\cb)$ is exact (whereas
$\ca\ten\cb$ is not, in general).

If $\ca$ is an arbitrary small dg category, there is a universal
dg functor
\[
\ca\to \pretr(\ca)
\]
to a pretriangulated dg category $\pretr(\ca)$, \ie a functor
inducing an equivalence
\[
\HOM(\ca,\cb) \iso \HOM(\pretr(\ca),\cb)
\]
for each exact dg category $\cb$. The dg category $\pretr(\ca)$ is
the {\em pretriangulated hull of $\ca$} constructed explicitly in
\cite{BondalKapranov90}, \cf also \cite{Drinfeld04}
\cite{Tabuada05b}.

For any dg category $\ca$, the category $H^0(\pretr(\ca))$ is
equivalent to the triangulated subcategory of $\ch\ca$ generated
by the representable dg modules. The functor $\pretr$ preserves
quasi-equivalences and induces a left adjoint to the inclusion of
the full subcategory of exact dg categories into the homotopy
category $\Hqe$. If $\cb$ is pretriangulated, then so is
$\RHOM(\ca,\cb)$ for each small dg category $\ca$ and we have
\[
\RHOM(\pretr(\ca),\cb) \iso \RHOM(\ca,\cb).
\]

\subsection{Morita fibrant dg categories, exact sequences}
\label{ss:MoritaFibrant} A dg functor $F:\ca \to \cb$ between small
dg categories is a {\em Morita morphism} if it induces an
equivalence $\cd(\cb)\to\cd(\ca)$. Each quasi-equivalence is a
Morita morphism (\cf section~\ref{ss:Morita-equivalence}) and so
is the canonical morphism $\ca\to \pretr(\ca)$ from $\ca$ to its
pretriangulated hull.

\begin{theorem}[\cite{Tabuada05b}] \label{thm:Tabuada2}
The category $\dgcat_k$ admits a structure of cofibrantly
generated model category whose weak equivalences are the Morita
morphisms and whose cofibrations are the same as those of the
canonical model structure on $\dgcat_k$ (\cf
theorem~\ref{thm:Tabuada1}).
\end{theorem}

A dg category $\ca$ is {\em Morita fibrant} (or {\em triangulated}
in the terminology of \cite{ToenVaquie05}) iff it is fibrant with
respect to this model structure. This is the case iff the
canonical functor $H^0(\ca) \to \per(\ca)$ is an equivalence iff
$\ca$ is pretriangulated and $H^0(\ca)$ is idempotent complete
(\cf section~\ref{ss:CompactObjectsBrownRep}). We write
$\ca\to\per_{dg}(\ca)$ for a fibrant replacement of $\ca$ and then
have
\[
\per(\ca) = H^0(\per_{dg}(\ca)).
\]
We write {\em $\Hmo$ for the localization} of $\dgcat_k$ with
respect to the Morita morphisms. Then the functor $\ca \mapsto
\per_{dg}(\ca)$ yields a right adjoint of the quotient functor
$\Hqe \to \Hmo$ and induces an equivalence from $\Hmo$ onto the
subcategory of Morita fibrant dg categories in $\Hqe$, \cf
\cite{Tabuada05b}. The category $\Hmo$ is pointed: The dg category
with one object and one morphism is both initial and terminal.
Moreover, $\Hmo$ admits all finite coproducts (they are induced by
the disjoint unions) and these are isomorphic to products.

Let
\begin{equation} \label{eq:ex-seq}
\xymatrix{\ca \ar[r]^I & \cb \ar[r]^{P} & \cc}
\end{equation}
be a sequence  of $\Hqe$ such that $PI=0$ in $\Hmo$.

\begin{theorem} The following are equivalent
\begin{itemize}
\item[i)] In $\Hmo$, $I$ is a kernel of $P$ and $P$ a cokernel of $I$.
\item[ii)] The morphism $I$ induces an equivalence of $\per(\ca)$
onto a thick subcategory of $\per(\cb)$ and $P$ induces an
equivalence of the idempotent closure \cite{BalmerSchlichting01}
of the Verdier quotient with $\per(\cc)$.
\item[iii)] The functor $I$ induces an equivalence of $\cd(\ca)$
with a thick subcategory of $\cd(\cb)$ and $P$ identifies the
Verdier quotient with $\cd(\cc)$.
\end{itemize}
\end{theorem}

The theorem is proved in \cite{Keller99}. The equivalence of ii)
and iii) is a consequence of Thomason-Trobaugh's localization
theorem \cite{ThomasonTrobaugh90} \cite{Neeman92a}
\cite{Neeman99}. We say that (\ref{eq:ex-seq}) is an {\em exact
sequence} of $\Hmo$ if the conditions of the theorem hold. For
example, if $X$ is a quasi-compact quasi-separated scheme,
$U\subset X$ a quasi-compact open subscheme and $Z=X\setminus U$,
then the sequence
\[
\parf_{dg}(X\  \mbox{on}\  Z) \to \parf_{dg}(X) \to \parf_{dg}(U)
\]
is an exact sequence of $\Hmo$ by the results of
\cite[Sect.~5]{ThomasonTrobaugh90}, where $\parf_{dg}(X)$ denotes
the dg quotient of the category of perfect complexes (viewed as a
full dg subcategory of the category of complexes of
$\co_X$-modules) by its subcategory of acyclic perfect complexes
and $\parf_{dg}(X\ \mbox{on}\ Z)$ the full subcategory of perfect
complexes supported on $Z$.

\subsection{Dg categories of finite type}
Let $\cm$ be a cofibrantly generated model category and $I$ a
small category. Recall that the category of functors $\cm^I$ is
again a cofibrantly generated model category (with the
componentwise weak equivalences). Thus, the diagonal functor
$\Ho(\cm)\to \Ho(\cm^I)$ admits a left adjoint, the {\em homotopy
colimit functor}, and a right adjoint, the {\em homotopy limit
functor}. An object $X$ of $\cm$ is {\em homotopically finitely
presented} if, for each filtered direct system $Y_i$, $i\in I$, of
$\cm$, the canonical morphism
\[
\hocolim \Map(X,Y_i) \to \Map(X, \hocolim Y_i)
\]
is a weak equivalence of simplicial sets. The category $\cm$ is
{\em homotopically locally finitely presented} if, in $\Ho(\cm)$,
each object is the homotopy colimit of a filtered direct system
(in $\cm$) of homotopically finitely presented objects.

For example \cite{Hinich97}, the category of dg algebras is
homotopically locally finitely presented and a dg algebra is
homotopically finitely presented iff, in the homotopy category, it
is a retract of a non commutative free graded algebra $k\langle
x_1, \ldots, x_n\rangle$ endowed with a differential such that
$dx_i$ belongs to $k\langle x_1, \ldots, x_{i-1}\rangle$ for each
$1\leq i\leq n$. A dg category is {\em of finite type} if it is dg
Morita equivalent to a homotopically finitely presented dg
algebra.

\begin{theorem}[\cite{ToenVaquie05}] \label{thm:finite-type}
The category of small dg categories endowed with the canonical
model structure whose weak equivalences are the Morita morphisms
is homotopically locally finitely presented and a dg category is
homotopically finitely presented iff it is of finite type.
\end{theorem}

A dg category $\ca$ is {\em smooth} if the bimodule $ (X,Y)
\mapsto \ca(X,Y) $ is perfect in $\cd(\ca\op\lten\ca)$. This
property is invariant under dg Morita equivalence. The explicit
description of the homotopically finitely presented dg algebras
shows that a dg category of finite type is smooth. Conversely
\cite{ToenVaquie05}, a dg category $\ca$ is of finite type if it
is smooth and {\em proper}, \ie dg Morita equivalent to a dg
algebra whose underlying complex of $k$-modules is perfect.

\subsection{Moduli of objects in dg categories} Let $T$ be a small dg
category. In \cite{ToenVaquie05}, B.~To\"en and M.~Vaqui\'e
introduce and study the $D^-$-stack (in the sense of
\cite{ToenVezzosi04}) of objects in $T$. By definition, this
$D^-$-stack is the functor
\[
\cm_T: \opname{Scalg}\to \Sset
\]
which sends a simplicial commutative $k$-algebra $A$ to the
simplicial set
\[
\Map(T\op, \per_{dg}(NA)) \ko
\]
where $NA$ is the commutative dg $k$-algebra obtained from $A$ by
the Dold-Kan equivalence. They show that if $T$ is a dg category
of finite type, then this $D^-$-stack is locally geometric and
locally of finite presentation. Moreover, if $E:T \to
\per_{dg}(k)$ is a $k$-point of $\cm_T$, then the tangent complex
of $\cm_T$ at $E$ is given by
\[
{\mathcal T}_{\cm_T, E} \iso \RHOM(E,E)[1].
\]
In particular, if $E$ is quasi-isomorphic to a representable
$x^\wedge$, then we have
\[
{\mathcal T}_{\cm_T, E} \iso T(x,x)[1].
\]
It follows that the restriction of $\cm_T$ to the category of
commutative $k$-algebras is a locally geometric $\infty$-stack in
the sense of C.~Simpson \cite{Simpson96}. Here are three
consequences derived from these results in \cite{ToenVaquie05}:

1) If $T$ is a dg category over a field $k$ and is smooth, proper
and Morita fibrant, then the sheaf associated with the presheaf
\[
R \mapsto \Aut_{\Hqe_R}(T\ten_k R) \ko
\]
on the category of commutative $k$-algebras is a group scheme
locally of finite type over $k$ (\cf \cite{Yekutieli04} for the
case where $T$ is an algebra).

2) If $X$ is a smooth proper scheme over a commutative ring $k$,
then the $\infty$-stack of perfect complexes on $X$ is locally
geometric.

3) If $A$ is a (non commutative) $k$-algebra over a field $k$,
then the $\infty$-stack of bounded complexes of finite-dimensional
$A$-modules is locally geometric if either $A$ is the path algebra
of a finite quiver or a finite-dimensional algebra of finite
global dimension.

\subsection{Dg orbit categories} Let $\ca$ be a dg category and
$F: \ca\to\ca$ an automorphism of $\ca$ in $\Hqe$. Let us assume
for simplicity that $F$ is given by a dg functor $\ca\to\ca$. The
{\em dg orbit category $\ca/F^\Z$} has the same objects as $\ca$
and the morphisms defined by
\[
(\ca/F^\Z)(X,Y) = \bigoplus_{d\in\Z}\colim_n \ca(F^{n}X,
F^{n+d}Y).
\]
The projection functor $P: \ca \to \ca/F^\Z$ is endowed with a
canonical morphism $\phi:PF \to P$ which becomes invertible in
$H^0(\ca/F^\Z)$ and the pair $(P,\phi)$ is the solution of a
universal problem, \cf \cite{Keller05}. The category
$H^0(\ca)/F^\Z$ is defined analogously. It is isomorphic to
$H^0(\ca/F^\Z)$ and can be thought of as the `category of orbits'
of the functor $F$ acting in $H^0(\ca)$.

Let us now assume that $k$ is a field. Let $Q$ be a quiver
(=oriented graph) whose underlying graph is a Dynkin graph of type
$A$, $D$ or $E$. Let $\mod kQ$ be the abelian category of
finite-dimensional representations of $Q$ over $k$ (\cf \eg
\cite{GabrielRoiter92} \cite{AuslanderReitenSmaloe95}). Let
$\ca=\cd^b_{dg}(\mod Q)$ and $F:\ca\to \ca$ an automorphism in
$\Hqe$. We say that $F$ {\em acts properly} if no indecomposable
object of $\cd^b(\mod kQ)$ is isomorphic to its image under $F$.
For example, if $\Sigma$ is the {\em Serre functor} of $\ca$,
defined by the bimodule
\[
(X,Y) \mapsto \HOM_k(\ca(Y,X),k)  \ko
\]
then $\Sigma$ acts properly and, more generally, if $S$ is the
suspension functor, then $S^{-d}\Sigma$ acts properly for each
$d\in\N$ unless $Q$ is reduced to a point.

\begin{theorem}[\cite{Keller05}]
If $F$ acts properly, the orbit category $\cd^b_{dg}(\mod
kQ)/F^\Z$ is Morita fibrant and thus $\cd^b(\mod kQ)/F^\Z$ is
canonically triangulated.
\end{theorem}

In the particular case where $F=S^{-d}\Sigma$, the triangulated
category $H^0(\ca/F^\Z)$ is Calabi-Yau \cite{Kontsevich98} of
CY-dimension $d$ (\cf \cite{Keller05}). For $d=1$, the category
$H^0(\ca/F^\Z)$ is equivalent to the category of
finite-dimensional projective modules over the preprojective
algebra (\cf \cite{GelfandPonomarev79} \cite{DlabRingel80}
\cite{Ringel98}) associated with the Dynkin graph underlying $Q$.
For $d=2$, one obtains the {\em cluster category} associated with
the Dynkin graph. This category was introduced in
\cite{CalderoChapotonSchiffler04} for type $A$ and in
\cite{BuanMarshReinekeReitenTodorov04} in the general case. It
serves in the representation-theoretic approach (\cf \eg
\cite{BuanMarshReinekeReitenTodorov04} \cite{CalderoChapoton04}
\cite{GeissLeclercSchroeer05}) to the study of cluster algebras
\cite{FominZelevinsky02} \cite{FominZelevinsky03}
\cite{BerensteinFominZelevinsky05} \cite{FominZelevinsky03a}. It
seems likely \cite{Amiot06} that if $k$ is algebraically closed,
the theorem yields `almost all' Morita fibrant dg categories whose
associated triangulated categories have finite-dimensional
morphism spaces and only finitely many isoclasses of
indecomposables. In particular, those among these categories which
are Calabi-Yau of fixed CY-dimension $d\gg 0$ are expected to be
parametrized by the simply laced Dynkin diagrams.

\section{Invariants}
\label{s:invariants}

\subsection{Additive invariants} Let $\Hmo_0$ be the category with
the same objects as $\Hmo$ and where morphisms $\ca\to\cb$ are
given by elements of the Grothendieck group of the triangulated
category $\rep(\ca,\cb)$. The composition is induced from that of
$\Hmo$. The category $\Hmo_0$ is additive and endowed with a
canonical functor $\Hmo\to\Hmo_0$ (\cf \cite{BondalLarsenLunts04}
for a related construction). One can show \cite{Tabuada05b} that a
functor $F$ defined on $\Hmo$ with values in an additive category
factors through $\Hmo\to\Hmo_0$ iff for each exact dg category
$\ca$ endowed with full exact dg subcategories $\cb$ and $\cc$
which give rise to a semi-orthogonal decomposition
$H^0(\ca)=(H^0(\cb),H^0(\cc))$ in the sense of
\cite{BondalKapranov90}, the inclusions induce an isomorphism
$F(\cb)\oplus F(\cc)\iso F(\ca)$. We then say that {\em $F$ is an
additive invariant}. The most basic additive invariant is given by
$F\ca=K_0(\per\ca))$. In $\Hmo_0$, it becomes a corepresentable
functor: $K_0(\per(\ca))= \Hmo_0(k,\ca)$. As we will see below,
the $K$-theory spectrum and all variants of cyclic homology are
additive invariants. This is of interest since non isomorphic
objects of $\Hmo$ can become isomorphic in $\Hmo_0$. For example,
if $k$ is an algebraically closed field, each finite-dimensional
algebra of finite global dimension becomes isomorphic to a product
of copies of $k$ in $\Hmo_0$ (\cf \cite{Keller98}) but it is
isomorphic to such a product in $\Hmo$ only if it is semi-simple.

\subsection{$K$-theory} \label{ss:K-theory}
Let $\ca$ be a small dg $k$-category. Its $K$-theory $K(\ca)$ is
defined by applying Waldhausen's construction \cite{Waldhausen85}
to a suitable category with cofibrations and weak equivalences:
here, the category is that of perfect $\ca$-modules, the
cofibrations are the morphisms $i: L \to M$ of $\ca$-modules which
admit retractions as morphisms of graded $\ca$-modules (\ie the
inflations of section~\ref{ss:triangulated-structure}) and the
weak equivalences are the quasi-isomorphisms. This construction
can be improved so as to yield a functor $K$ from $\dgcat_k$ to
the homotopy category of spectra. As in \cite{ThomasonTrobaugh90},
from Waldhausen's results \cite{Waldhausen85} one then obtains the
following

\begin{theorem}
\begin{itemize}
\item[a)] \cite{DuggerShipley04} The map $\ca \mapsto K(\ca)$ yields a well-defined
functor on $\Hmo$.
\item[b)] Applied to the bounded dg-derived category
$\cd^b_{dg}(\ce)$ of an exact category $\ce$, the $K$-theory
defined above agrees with Quillen $K$-theory.
\item[c)] The functor $\ca\mapsto K(\ca)$ is an additive
invariant. Moreover, each short exact sequence $\ca\to\cb\to\cc$
of $\Hmo$ (\cf section~\ref{ss:MoritaFibrant}) yields a long exact
sequence
\[
\ldots \to K_i(\ca) \to K_i(\cb) \to K_i(\cc) \to \ldots \to
K_0(\cb) \to K_0(\ca).
\]
\end{itemize}
\end{theorem}

Part a) can be improved on: In fact, D.~Dugger and B.~Shipley show
in \cite{DuggerShipley04} that $K$-theory is even preserved under
{\em topological} Morita equivalence. Part c) can be improved on
by defining {\em negative $K$-groups} and showing that the exact
sequence extends indefinitely to the right. We refer to
\cite{Schlichting06} for the most recent results, which include
the case of dg categories. By combining part a) with Rickard's
theorem~\ref{thm:Rickard}, one obtains the invariance of the
$K$-theory of rings under triangle equivalences between their
derived categories. By combining a) and b), one obtains the
invariance of the $K$-theory of abelian categories under
equivalences between their derived categories which come from
isomorphisms of $\Hmo$ (or, more generally, from topological
Morita equivalences). In fact, according to A.~Neeman's results
\cite{Neeman97}
the $K$-theory of an abelian category is even determined by the
underlying triangulated category of its derived category, \cf
\cite{Neeman05} for a survey of his work.

Of course, any invariant defined for small triangulated categories
applied to the perfect derived category yields an invariant of
small dg categories. For example, Balmer-Witt groups (\cf
\cite{Balmer04} for a survey), defined for dg categories $\ca$
endowed with a suitable involution $\ca\iso\ca\op$ in $\Hmo$,
yield such invariants.

\subsection{Hochschild and cyclic homology}
Let $\ca$ be a small $k$-flat $k$-category. Following
\cite{Mitchell72} the {\em Hochschild chain complex} of $\ca$ is
the complex concentrated in homological degrees $p\geq 0$ whose
$p$th component is the sum of the
\[
\ca(X_p,X_0)\ten \ca(X_p,X_{p-1}) \ten \ca(X_{p-1}, X_{p-2}) \ten
\cdots \ten \ca(X_0,X_1)\ko
\]
where $X_0, \ldots, X_p$ range through the objects of $\ca$,
endowed with the differential
\[
d(f_{p}\ten \ldots \ten f_0) = f_{p-1}\ten\cdots \ten f_0 f_p +
\sum_{i=1}^p (-1)^i f_p\ten \cdots \ten f_i f_{i-1} \ten\cdots\ten
f_0.
\]
Via the cyclic permutations
\[
t_p (f_{p-1} \ten \cdots \ten f_0) = (-1)^p f_0 \ten f_{p-1} \ten
\cdots \ten f_1
\]
this complex becomes a precyclic chain complex and thus gives rise
\cite[Sect.~2]{Keller96} to a {\em mixed complex $C(\ca)$} in the
sense of \cite{Kassel85}, \ie a dg module over the dg algebra
$\Lambda=k[B]/(B^2)$, where $B$ is of degree $-1$ and $dB=0$. As
shown in \cite{Kassel85}, all variants of cyclic homology
\cite{Loday98} only depend on $C(\ca)$ considered in
$\cd(\Lambda)$. For example, the cyclic homology of $\ca$ is the
homology of the complex $C(\ca)\lten_\Lambda k$.

If $\ca$ is a $k$-flat differential graded category, its mixed
complex is the sum-total complex of the bicomplex obtained as the
natural re-interpretation of the above complex. If $\ca$ is an
arbitrary dg $k$-category, its Hochschild chain complex is defined
as the one of a $k$-flat (\eg a cofibrant) resolution of $\ca$.

\begin{theorem}[\cite{Keller98a} \cite{Keller99}]
\begin{itemize}
\item[a)] The map $\ca \mapsto C(\ca)$
yields an additive functor $\Hmo_0 \to \cd(\Lambda)$. Moreover,
each exact sequence of $\Hmo$ (\cf section~\ref{ss:MoritaFibrant})
yields a canonical triangle of $\cd(\Lambda)$.
\item[b)] If $A$ is a $k$-algebra, there is a natural
isomorphism $C(A) \iso C(\per_{dg}(A))$.
\item[c)] If $X$ is a quasi-compact separated scheme, there is a
natural isomorphism $C(X) \iso C(\parf_{dg}(X))$, where $C(X)$ is
the cyclic homology of $X$ in the sense of \cite{Loday86}
\cite{Weibel96} and $\parf_{dg}(X)$ the dg category defined in
section~\ref{ss:MoritaFibrant}.
\end{itemize}
\end{theorem}

The second statement in a) may be viewed as an excision theorem
analogous to \cite{Wodzicki89}. We refer to the recent proof
\cite{CortinasHaesemeyerSchlichtingWeibel05} of Weibel's
conjecture \cite{Weibel80} on the vanishing of negative $K$-theory
for an application of the theorem. The algebraic description of
{\em topological} Hochschild (co-)homology \cite{Shipley00} would
suggest that it is also preserved under topological Morita
equivalence but no reference seems to exist as yet.

The endomorphism algebra $\RHOM_\Lambda(k,k)$ is quasi-isomorphic
to $k[u]$, where $u$ is of degree $2$ and $d(u)=0$. It acts on
$C(\ca)\lten_{\Lambda} k$ and this action is made visible in the
isomorphism
\[
C(\ca)\lten_{\Lambda} k = C(\ca)\ten k[u]
\]
where $u$ is of degree $2$ and the differential on the right hand
complex is given by
\[
d(x\ten f) = d(x)\ten f + (-1)^{|x|}xB \ten uf.
\]
The following `Hodge--de Rham conjecture' is true for the dg
category of perfect complexes on a smooth projective variety or
over a finite-dimensional algebra of finite global dimension. It
is wide open in the general case.

\begin{conjecture}[\cite{Drinfeld02a} \cite{Kontsevich04}]
If $\ca$ is a smooth proper dg category over a field $k$ of
characteristic $0$, then the homology of $C(\ca)\ten k[u]/(u^n)$
is a flat $k[u]/(u^n)$-module for all $n\geq 1$.
\end{conjecture}

\subsection{Hochschild cohomology}
\label{ss:Hochschild-cohomology} Let $\ca$ be a small cofibrant dg
category. Its cohomological Hochschild complex $C(\ca,\ca)$ is
defined as the product-total complex of the bicomplex whose $0$th
column is
\[
\prod \ca(X_0,X_0) \ko
\]
where $X_0$ ranges over the objects of $\ca$, and whose $p$th
column, for $p\geq 1$, is
\[
\prod \HOM_k(\ca(X_{p-1}, X_p)\ten \ca(X_{p-2},X_{p-1})
\ten\cdots\ten \ca(X_{0},X_1), \ca(X_0,X_p))
\]
where $X_0, \ldots, X_p$ range over the objects of $\ca$. The
horizontal differential is given by the Hochschild differential.
This complex carries rich additional structure: As shown in
\cite{GetzlerJones94}, it is a $B_\infty$-algebra, \ie its bar
construction carries, in addition to its canonical differential
and comultiplication, a natural {\em multiplication} which makes
it into a dg bialgebra. The $B_\infty$-structure contains in
particular the cup product and the Gerstenhaber bracket, which
both descend to the Hochschild cohomology
\[
\HHs(\ca,\ca)=H^* C(\ca,\ca).
\]
The Hochschild cohomology is naturally interpreted as the homology
of the complex
\[
\HOM(\id_\ca, \id_\ca)
\]
computed in the dg category $\RHOM(\ca,\ca)$, where $\id_\ca$
denotes the identity functor of $\ca$ (\ie the bimodule
$(X,Y)\mapsto \ca(X,Y)$). Then the cup product corresponds to the
composition (whereas the Gerstenhaber bracket has no obvious
interpretation). Each $c \in HH^n(\ca,\ca)$ gives rise to
morphisms $cM: M \to M[n]$ of $\cd(\ca)$, functorial in $M\in
\cd(\ca)$. Another interpretation links the Hochschild cohomology
of $\ca$ to the derived Picard group and to the higher homotopy
groups of the category of quasi-equivalences between dg
categories, \cf section~\ref{ss:InvertingQuasiEq}.

A natural way of obtaining the $B_\infty$-algebra structure on
$C(\ca,\ca)$ is to consider the $A_\infty$-category of
$A_\infty$-functors from $\ca$ to itself \cite{Kontsevich98}
\cite{Lefevre03} \cite{Lyubashenko03}. Here, the
$B_\infty$-algebra $C(\ca,\ca)$ appears as the endomorphism
algebra of the identity functor (\cf \cite{Keller05c}).

Note that $C(\ca,\ca)$ is not functorial with respect to dg
functors. However, if $F:\ca\to\cb$ is a fully faithful dg
functor, it clearly induces a restriction map
\[
F^* : C(\cb,\cb) \to C(\ca,\ca)
\]
and this map is compatible with the $B_\infty$-structure. This can
be used to construct \cite{Keller03} a morphism
\[
\phi_X : C(\cb,\cb) \to C(\ca,\ca)
\]
in the homotopy category of $B_\infty$-algebras associated with
each dg $\ca$-$\cb$-bimodule $X$ such that the functor
\[
?\lten_\ca X : \per(\ca) \to \cd\cb
\]
is fully faithful. If moreover the functor $ X\lten_\cb ? :
\per(\cb\op) \to \cd(\ca\op) $ is fully faithful, then $\phi_X$ is
an isomorphism. In particular, the Hochschild complex becomes a
functor
\[
\Hmo_{\mbox{\scriptsize ff}}\op \to \Ho(B_\infty) \ko
\]
where $\Ho(B_\infty)$ is the homotopy category of
$B_\infty$-algebras and $\Hmo_{\mbox{\scriptsize ff}}$ the (non
full) subcategory of $\Hmo$ whose morphisms are the quasi-functors
$X\in \rep(\ca,\cb)$ such that
\[
?\lten_\ca X : \per(\ca) \to \per(\cb)
\]
is fully faithful. We refer to \cite{LowenVandenBergh04a} for the
closely related study of the Hochschild complex of an abelian
category.

Let us suppose that $k$ is a field of characteristic $0$. Endowed
with the Gerstenhaber bracket the Hochschild complex $C(\ca,\ca)$
becomes a differential graded Lie algebra and this Lie algebra
`controls the deformations of the $A_\infty$-category $\ca$', \cf
\eg \cite{KontsevichSoibelman01}. Here the $A_\infty$-structures
$(m_n)$, $n\geq 0$, may have a non trivial term $m_0$. Some (but
not all) Hochschild cocyles also correspond to deformations of
$\ca$ as an object of $\Hmo$. To be precise, let $k[\eps]$ be the
algebra of dual numbers and consider the reduction functor
\[
R : \Hmo_{k[\eps]} \to \Hmo_k \ko \cb \mapsto \cb \lten_{k[\eps]}
k .
\]
A {\em first order Morita deformation} of $\ca$ is a pair $(\ca',
\phi)$ formed by a dg $k[\eps]$-category $\ca'$ and an isomorphism
$\phi: R\ca' \to \ca$ of $\Hmo_k$. An equivalence between such
deformations is given by an isomorphism $\psi: \ca'\to\ca''$ such
that $\phi' R\psi = \phi$. Then one can show \cite{GeissKeller05a}
that the equivalence classes of first order Morita deformations of
$\ca$ are in natural bijection with the classes $c \in
HH^2(\ca,\ca)$ such that the induced morphism $ cP : P \to P[2] $
is nilpotent in $H^*\HOM(P,P)$ for each perfect $\ca$-module $P$.
If $\ca$ is proper or, more generally, if $H^n\ca(?,X)$ vanishes
for $n\gg 0$ for all objects $X$ of $\ca$, then this condition
holds for all Hochschild $2$-cocycles $c$. On the other hand, if
$\ca$ is given by the dg algebra $k[u,u^{-1}]$, where $u$ is of
degree $2$ and $du=0$, then it does not hold for the cocycle $u\in
HH^2(\ca,\ca)$.

\subsection{Fine structure of the Hochschild complexes}
The Hochschild cochain complex of a dg category carries a natural
homotopy action of the little squares operad. This is the positive
answer to a question by P.~Deligne \cite{Deligne93} which has been
obtained, for example, in \cite{McClureSmith02}
\cite{KontsevichSoibelman00} \cite{BergerFresse04} \ldots.
Hochschild cohomology acts on Hochschild homology and this action
comes from a homotopy action of the Hochschild cochain complex,
viewed as a homotopy algebra over the little squares, on the
Hochschild chain complex. This is the positive answer to a series
of conjectures due to B.~Tsygan \cite{Tsygan99} and
Tamarkin-Tsygan \cite{TamarkinTsygan00}. It has recently been
obtained by B.~Tsygan and D.~Tamarkin \cite{Tsygan04}. Together,
the two Hochschild complexes endowed with these structures yield a
{\em non commutative calculus} \cite{TamarkinTsygan05} analogous
to the differential calculus on a smooth manifold. The link with
classical calculus on smooth commutative manifolds is established
through M.~Kontsevich's formality theorem \cite{Kontsevich97}
\cite{Tamarkin98} for Hochschild cochains and in \cite{Shoikhet03}
(\cf also \cite{Dolgushev03}) for Hochschild chains.

Clearly, these finer structures on the Hochschild complexes are
linked to the category of dg categories and its simplicial
enrichment given by the Dwyer-Kan localization as developped in
\cite{Toen04}. At the end of the introduction to \cite{Toen04},
the reader will find a more detailed discussion of these links,
\cf also \cite{KockToen05}. A precise relationship is announced in
\cite{Tamarkin05}.

\bigskip

\subsection{Derived Hall algebras} \label{ss:Hall-algebras}
Let $\ca$ be a {\em finitary} abelian category, \ie such that the
underlying sets of $\ca(X,Y)$ and $\Ext^1(X,Y)$ are finite for all
objects $X$, $Y$ of $\ca$. The {\em Ringel-Hall algebra
$\ch(\ca)$} is the free abelian group on the isomorphism classes
of $\ca$ endowed with the multiplication whose structure constants
are given by the Hall numbers $f_{XY}^Z$, which count the number
of subobjects of $Z$ isomorphic to $X$ and such that $Z/X$ is
isomorphic to $Y$, \cf \cite{DengXiao04} for a survey. Thanks to
Ringel's famous theorem \cite{Ringel90} \cite{Ringel93}, for each
simply laced Dynkin diagram $\Delta$, the {\em positive part} of
the Drinfeld-Jimbo quantum group $U_q(\Delta)$ (\cf \eg
\cite{Kassel95} \cite{Lusztig93}) is obtained as the (generic, twisted) Ringel-Hall
algebra of the abelian category of finite-dimensional
representations of a quiver $\vec{\Delta}$ with underlying graph
$\Delta$. After Ringel's discovery, it was first pointed out by
Xiao \cite{Xiao95}, \cf also \cite{Kapranov98}, that an extension of
the construction of the Ringel-Hall algebra to the derived
category of the representations of $\vec{\Delta}$ might yield the
{\em whole} quantum group. However, if one tries to mimick the
construction of $\ch(\ca)$ for a triangulated category $\ct$ by
replacing short exact sequences by triangles one obtains a
multiplication which fails to be associative, \cf\
\cite{Kapranov98} \cite{Hubery04}. It is remarkable
that nevertheless, as shown by Peng-Xiao \cite{PengXiao96}
\cite{PengXiao97} \cite{PengXiao00}, the
commutator associated with this multiplication yields the correct Lie
algebra.

A solution to the problem of constructing an associative
multiplication from the triangles has recently been proposed by B.~To\"en
in \cite{Toen05}.  He obtains an explicit formula for the
structure constants $\phi_{XY}^Z$ of an associative multiplication
on the rational vector space generated by the isomorphism classes
of any triangulated category $\ct$ which appears as the perfect
derived category $\per(T)$ of a proper dg category $T$ over a
finite field $k$. The resulting $\Q$-algebra is the {\em derived
Hall algebra $\derh(T)$ of $T$}. The formula for the structure
constants reads as follows:
\[
\phi_{XY}^Z= \sum_{f} |\Aut(f/Z)|^{-1} \prod_{i>0}
|\Ext^{-i}(X,Z)|^{(-1)^i} |\Ext^{-i}(X,X)|^{(-1)^{i+1}}\;\;,
\]
where $f$ ranges over the set of orbits of the group $\Aut(X)$ in
the set of morphisms $f: X \to Z$ whose cone is isomorphic to $Y$,
and $\Aut(f/Z)$ denotes the stabilizer of $f$ under the action of
$\Aut(X)$. The proof of associativity is inspired by  methods from
the study of higher moduli spaces \cite{ToenVezzosi04}
\cite{Toen04} \cite{ToenVaquie05} and by the homotopy theoretic
approach to $K$-theory \cite{Quillen73}. From the formula, it is
immediate that $\derh(T)$ is preserved under triangle equivalences
$\per(T) \iso \per(T')$. Another consequence is that if $\ca$ is
the heart of a non degenerate $t$-structure
\cite{BeilinsonBernsteinDeligne82} on $\per(T)$, then the
Ringel-Hall algebra of $\ca$ appears as a subalgebra of
$\derh(T)$. The derived Hall algebra of the derived category of
representations of $\vec{\Delta}$ over a finite field appears
closely related to the constructions of \cite{Kapranov98}. Its
precise relation to the quantum group $U_q(\Delta)$ remains to be
investigated.

Notice that like the $K_0$-group, the derived Hall algebra only
depends on the underlying triangulated category of $\per(T)$. 
One would expect that geometric versions of the
derived Hall algebra, as defined in \cite[3.3]{Toen05a}
will depend on finer data.


\def\cprime{$'$}
\providecommand{\bysame}{\leavevmode\hbox to3em{\hrulefill}\thinspace}
\providecommand{\MR}{\relax\ifhmode\unskip\space\fi MR }
\providecommand{\MRhref}[2]{%
  \href{http://www.ams.org/mathscinet-getitem?mr=#1}{#2}
}
\providecommand{\href}[2]{#2}

\end{document}